\documentclass[10pt]{article}
\usepackage{amsmath,amsthm,amssymb,amscd,fancybox,ifthen,float,epsfig,subcaption}

\usepackage[absolute,overlay]{textpos}

\usepackage{ragged2e}
\usepackage{tabto}

\floatstyle{plain}
\restylefloat{figure}

\usepackage{multirow}
\usepackage{tikz}
\usepackage{color}
\usepackage{url}
\usepackage{afterpage}
\usepackage{multicol}

\usepackage{marginnote}
\usepackage{tcolorbox}

\usepackage{lmodern}

\usepackage[sf,bf,medium]{titlesec}

\usepackage{bm}
\usepackage{booktabs}
\usepackage{siunitx}
\usepackage{isomath}

\usepackage[plainpages=false, colorlinks=true, citecolor=blue, filecolor=blue,
   linkcolor=blue, urlcolor=blue]{hyperref}
\usepackage{enumitem}

\usepackage{multicol}
\newcommand{\lp}{\left(}
\newcommand{\rp}{\right)}

\usepackage[sort,compress]{cite}

\newcommand\bbR{\mathbb R}
\newcommand\bbC{\mathbb C}
\newcommand\pa{\partial}

\newcommand\bx{\boldsymbol{x}}

\newcommand\by{\boldsymbol{y}}
\newcommand\bz{\boldsymbol{z}}

\newcommand\In{\operatorname{inc}}
\newcommand\inc{\operatorname{inc}}
\newcommand\PV{\operatorname{P.V.}}
\newcommand\out{\operatorname{out}}
\newcommand\rad{\operatorname{rad}}
\newcommand\tot{\operatorname{tot}}
\newcommand\loc{\operatorname{loc}}

\newcommand\inskel{\operatorname{in-skel}}
\newcommand\bn{\boldsymbol n}

\newcommand\tL{\widetilde L}

\newtheorem{theorem}{\sffamily Theorem}
\newtheorem{remark}{\sffamily Remark}

\newtheorem{lemma}{\sffamily Lemma}

\newcommand{\cD}{\mathcal D}
\newcommand{\cS}{\mathcal S}

\newcommand{\cC}{\mathcal C}

\newcommand{\cG}{\mathcal G}

\newcommand{\fd}{\mathfrak{d}}
\newcommand{\fv}{\mathfrak{v}}

\newcommand{\Id}{\rm{Id}}

\numberwithin{equation}{section}
\newcommand{\rhs}{\operatorname{rhs}}

\newcommand{\bs}{\boldsymbol}

\newcommand{\hsigma}{\hat\sigma}
\newcommand{\htau}{\hat\tau}

\newcommand{\figref}[1]{\figurename~\hyperref[#1]{\ref{#1}}}
\newcommand{\tableref}[1]{\tablename~\hyperref[#1]{\ref{#1}}}

\newcommand\tw{\widetilde w}
\newcommand\tu{\widetilde u}
\newcommand\tV{\widetilde V}
\newcommand\supp{\operatorname{supp}}

\renewcommand{\phi}{\varphi}

\newcommand{\qand}{\quad\text{and}\quad}

\binoppenalty=\maxdimen
\relpenalty=\maxdimen

\newcommand{\nchnk}{{n_{\rm{chunk}}}}

\usepackage{setspace}
\usetikzlibrary{calc}
\usepackage{abstract}

\begin{document}

\begin{titlepage}

  \raggedleft
  {\sffamily \bfseries STATUS: arXiv pre-print}

  \raggedright
 \begin{textblock*}{\linewidth}(1.25in,2in) 
    {\LARGE \sffamily \bfseries A numerical method for scattering problems with unbounded interfaces}
  \end{textblock*}

   \vspace{1.5in}
   Tristan Goodwill\\ \textit{\small Department of Statistics and CCAM\\
   University of Chicago\\
   Chicago, IL 60637}\\ \texttt{\small tgoodwill@uchicago.edu}


   \vspace{1.5\baselineskip}
    Charles L. Epstein
  \\ \textit{\small Center for Computational Mathematics\\
   Flatiron Institute, Simons Foundation\\
   New York, NY 10010}\\ \texttt{\small cepstein@flatironinstitute.org}

  \begin{textblock*}{\linewidth}(1.25in,7in) 
    \today
  \end{textblock*}

\end{titlepage}

\begin{abstract}
  \noindent 
  We introduce a new class of computationally tractable scattering problems in unbounded domains, which we call \emph{decomposable problems.} In these decomposable problems, the computational domain can be split into a finite collection of subdomains in which the scatterer has a ``simple" structure. A subdomain is simple if  the domain Green's function for this subdomain is either available analytically or can be computed numerically with arbitrary accuracy by a tractable method. These domain Green's functions are then used to reformulate the scattering problem as a system of boundary integral equations on the union of the subdomain boundaries. This reformulation  gives a practical numerical method, as the resulting integral equations can then be solved, to any desired degree of accuracy, by using  coordinate complexification over a \emph{finite} interval, and standard discretization techniques.\\


\end{abstract}

{
  \hypersetup{linkcolor=black}
  \tableofcontents
}


\normalsize

\section{Introduction}
\begin{figure}[h]
  \centering
      \includegraphics[width=8cm]{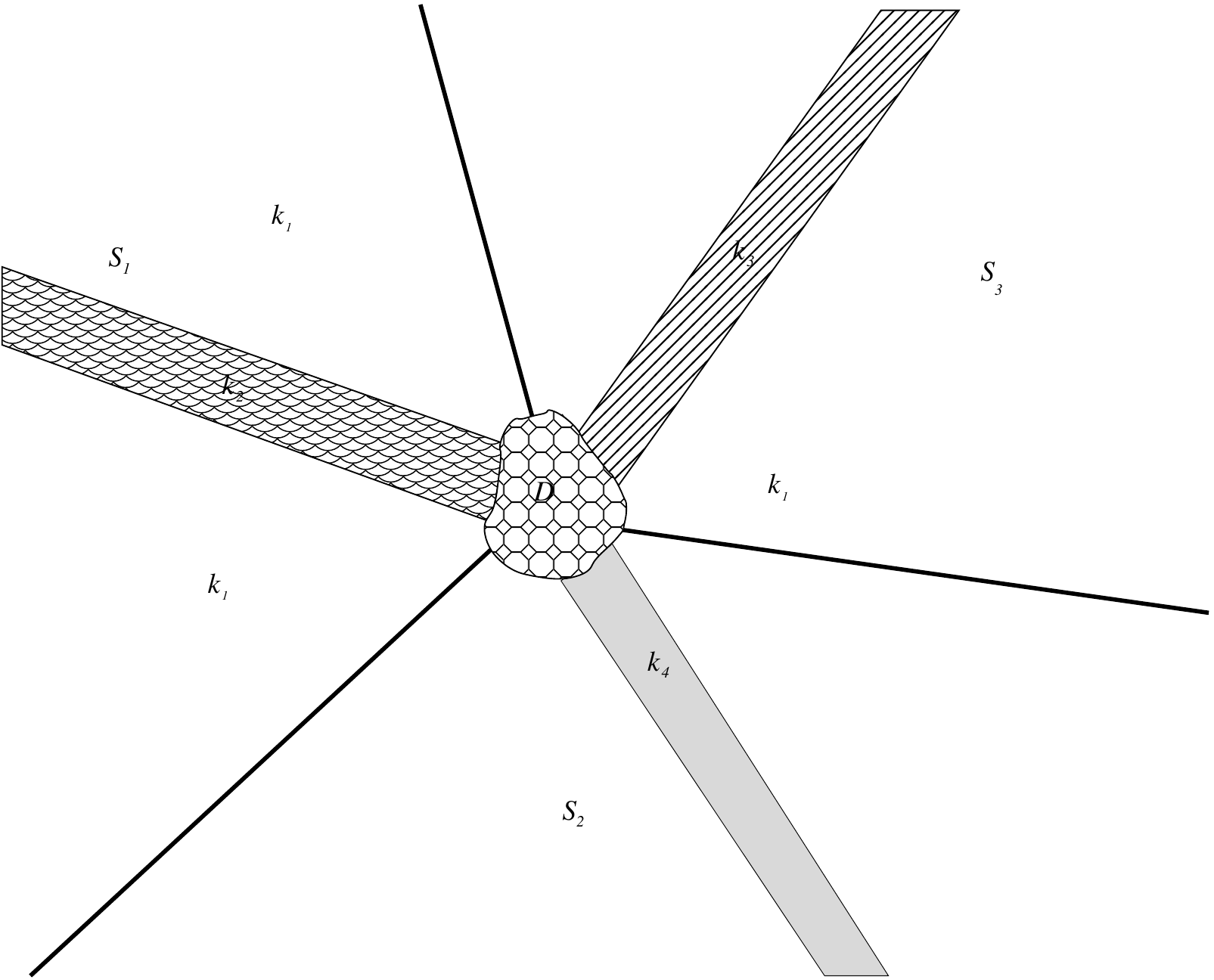}
      \caption{Three dielectric channels meeting in a compact interaction zone,
        $D,$ showing sectors $S_1, S_2,S_3.$}
       \label{fig1}
      \end{figure}

There are many circumstances where one would like to examine how waves scatter
off of unbounded interfaces. In this paper we consider only scalar waves
described as perturbations of the Helmholtz equation.  A
typical unbounded scatterer is a network of  open wave-guides described by
$\Delta+q({\bx})+k^2$ where the
potential is of the form
\begin{equation}\label{eqn1}
  q({\bx})=q_0({\bx})+\sum_{j=1}^Nq_j({\bx}-\langle {\bx},\omega_j\rangle \omega_j)
  \chi_{[1,\infty)}(r_j\langle {\bx},\omega_j\rangle).
\end{equation}
Here $q_0$ is a compactly supported function, the $\{q_j:\:j=1,N\}$ are
compactly supported functions of one fewer variable, the $\{\omega_j\}$ are distinct
unit vectors, and the $\{r_j\}$ are positive numbers. A schematic example is shown in Figure~\ref{fig1}.

\begin{figure}[h]
  \centering
    \includegraphics[width= 8cm]{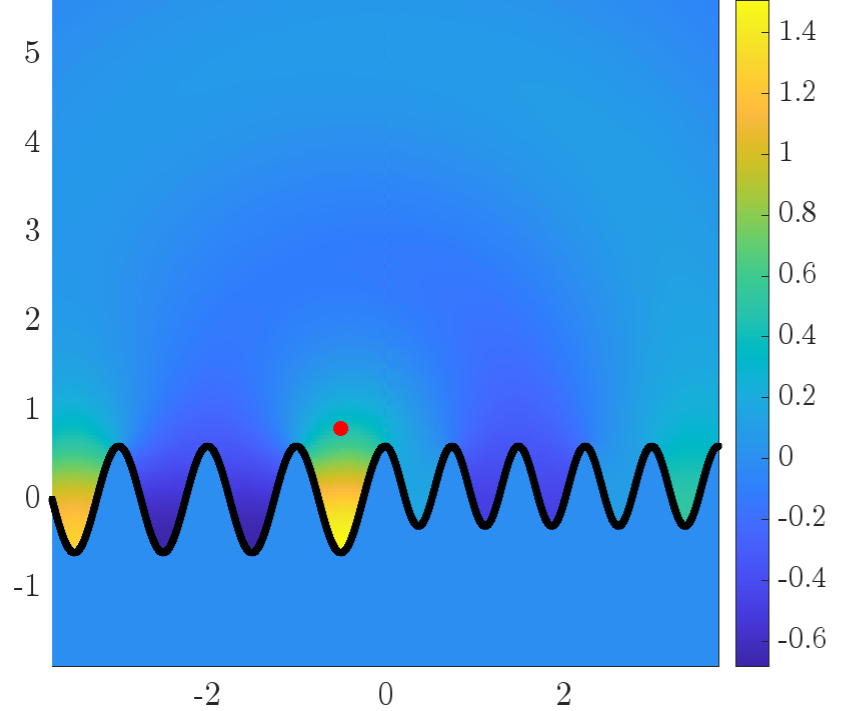}
    \caption{Two semi-infinite periodic boundaries meeting at along a common
      perpendicular line.}
    \label{fig2}
    \end{figure}

Another situation where
such interfaces arise is that of a half space with either a potential supported
near the boundary, or perhaps a boundary that is itself non-flat over a
non-compact set, see Figure~\ref{fig2}. Problems of this sort have proved quite resistant to accurate
numerical solution.  

Methods for this class of problem can be placed into one of a few categories. Many use approximations to truncate the domain~\cite{ hadley1992transparent, huang1993simulation, kim2015optimized, koshiba2001high, orlandini2024waveguide,zhang2011novel,LuLuQian,bourgeois2023scattering,shibayama2006performance,becache2021stability}. Of particular interest is the method described in~\cite{bruno2017windowed}, which converge superalgebraically with truncation length. There also exist methods based on time domain calculations~\cite{sofronov2015application,alpert2002nonreflecting,mulder2020exact}, or are restricted to the case where~\eqref{eqn1} is a perturbation of a layered media problem~\cite{cai2013computational,okhmatovski2004evaluation,wang2019fast,zhang2020exponential,paulus2000accurate,taraldsen2005complex,perez2014high,bruno2016windowed}. There also exist a few methods that have solved the problem by imposing artificial radiation conditions~\cite{BonnetBendhia_etal,ott2017halfspace}. In this and subsequent papers in this series, we present a new approach to the numerical solution of a certain subclass of such problems. This method has been shown to give the unique solution of the problem with physically meaningful radiation conditions and gives provably exponentially small truncation error. Further, the reformulated method can easily be discretized using standard integral equation techniques.

We let $L=\Delta+q({\bx});$ the scattering
problem is to find an ``outgoing'' solution to
\begin{equation}\label{eqn1.2}
  (L+k^2)u=f,
\end{equation}
where $f$ is defined by ``incoming data.'' The precise formulation of incoming and outgoing
radiation conditions depends on the problem at hand. The operator $L$ is usually
self-adjoint so that the resolvent operator $(L+k^2+i\delta)^{-1}$ is a bounded
operator on $L^2$ provided that $\delta\neq 0.$ The limiting absorption
principle then states that the limits
\begin{equation}
  (L+k^2\pm i0)^{-1}=\lim_{\delta\to 0^{\pm}}(L+k^2+i\delta)^{-1}
\end{equation}
exist as bounded operators between weighted spaces, e.g. as a bounded maps
$$ (L+k^2\pm i0)^{-1}:(1+|{\bx}|^2)^{-\frac{1+\eta}{4}}L^2(\bbR^d)\longrightarrow
(1+|{\bx}|^2)^{\frac{1+\eta}{4}}L^2(\bbR^d),$$
for any $\eta>0,$  see~\cite{Vasy2000}. If $f$ is compactly supported, then $u^+=(L+k^2+ i0)^{-1}f$ is
typically the unique outgoing solution to $(L+k^2)u^+=f.$
The long-range nature of the potential and the existence of wave-guide modes
make it much more complicated to state the radiations conditions, and check that
they imply uniqueness. Physically motivated radiated conditions in
$d$-dimensions, implying uniqueness for potentials like those in~\eqref{eqn1},
are given in~\cite{EpMaSR2023}.

If $q({\bx})=(k_1^2-k^2)\chi_{D}({\bx}),$ for $D$ a compact region with a smooth
boundary, then the standard scattering problem is formulated as follows: We let
$v^{\In}$ be a solution to $(\Delta+k^2)v^{\In}=0$ in the exterior of $D,$ 
smooth up to $\pa D.$ We then look for functions $u_{\pm}$ with
\begin{equation}
  \begin{split}
    &(\Delta+k^2)u_+({\bx})=0\text{ for }{\bx}\in D^c,\\
    &(\Delta+k_1^2)u_-({\bx})=0\text{ for }{\bx}\in D,    
  \end{split}
\end{equation}
and set
\begin{equation}
  u^{\tot}({\bx})=
  \begin{cases}
    v^{\In}({\bx})+u_+({\bx})&\text{ for }{\bx}\in D^c,\\
    u_-({\bx})&\text{ for }{\bx}\in D.
  \end{cases}
\end{equation}
In addition, $u_+$ is required to satisfy the outgoing radiation condition,
which in 2-dimensions means that
\begin{equation}
  (\pa_r-ik)u_+(r\omega)=o(r^{-\frac 12}),
\end{equation}
uniformly in $\omega\in S^1.$ The solution is uniquely determined by the condition that
$u^{\tot}\in H^2_{\loc}(\bbR^2),$ which is equivalent to the transmission
boundary conditions
\begin{equation}\label{eqn1.7.02}
  \begin{split}
    u_-({\bx})-u_+({\bx})&=v^{\In}({\bx})\text{ for }{\bx}\in\pa D,\\
    \pa_{\bn}u_-({\bx})-\pa_{\bn}u_+({\bx})&=\pa_{\bn}v^{\In}({\bx})\text{ for }{\bx}\in\pa D,
  \end{split}
\end{equation}
here $\bn$ is the outer unit normal vector field to $\pa D.$

The standard method to solve this problem is to express $u_{\pm}$ in terms of
layer potentials on $\pa D.$ The function
\begin{equation}
  G_{k}({\bx};{\by})=\frac{i}{4\pi}H^{(1)}_0(k|{\bx}-{\by}|),
\end{equation}
is the outgoing fundamental solution for $(\Delta+k^2);$  for ${\bx}\notin \pa D,$ we let
\begin{equation}
  \begin{split}
    &\cS_k [\tau]({\bx})=\int_{\pa D}G_k({\bx};{\by})\tau({\by})ds({\by})\\
    &\cD_k [\sigma]({\bx})=\int_{\pa D}\pa_{\bn_{\by}}G_k({\bx};{\by})\sigma({\by})ds({\by})
    \end{split}
  \end{equation}
  denote the single and double layer potentials. Assume that
  \begin{equation}
    \begin{split}
    & u_+({\bx})=\cS_k[\tau]({\bx})+\cD_k[\sigma]({\bx})\text{ for }{\bx}\in D^c,\\
      & u_-({\bx})=\cS_{k_1}[\tau]({\bx})+\cD_{k_1}[\sigma]({\bx})\text{ for }{\bx}\in D.
      \end{split}
  \end{equation}
  Restricting $u_{\pm}$ to the $\pa D$ and using standard jump conditions for layer
  potentials we find that the conditions in~\eqref{eqn1.7.02} are equivalent to
  the system of integral equations
  \begin{equation}
    \left(  \begin{matrix}
      \Id + K_1& K_2\\K_3&\Id + K_4
    \end{matrix}
\right)\left(  \begin{matrix}
      \sigma\\
      \tau
    \end{matrix}
\right)=  \left(  \begin{matrix}
      -v^{\In}\\
      \pa_{\bn}v^{\In}
    \end{matrix}
\right),
  \end{equation}
  where the $\{K_j\}$ are compact operators acting on H\"older continuous
  functions on $\pa D.$ As such these are Fredholm equations of second kind, and
  therefore generally well conditioned.
  By construction, $u_+$ is an outgoing solution.
  
In this work, we adapt this method  to obtain numerical solutions of scattering problems from non-compact
scatterers. The principal idea  is to geometrically
decompose the problem into sub-problems for which the outgoing fundamental
solutions can be efficiently  approximated numerically. These
fundamental solutions are no longer available in closed form,  but require the solutions
of  non-trivial, but tractable, PDE problems. Cases where these sorts of more
complicated Green's
functions have been used include layered media with piecewise
constant layers~\cite{cai2013computational,okhmatovski2004evaluation,wang2019fast,zhang2020exponential,paulus2000accurate,taraldsen2005complex,bruno2016windowed}, and subsets of a halfspace with a periodic boundary~\cite{AgocsBarnett_2024}. 
 Our approach to such scattering problems has three principal components.
\begin{enumerate}
\item Replace the scattering problem with a transmission
  problem. 
  \begin{enumerate}
  \item Divide space into a finite collection of half spaces, cylinders, and asymptotically
  conical sectors $\{S_j:\: j=1,\dots,N\},$ along with a compact set $S_0$ (see Figure~\ref{fig1}). These
  subsets meet along piecewise smooth curves that lie on straight lines outside of a
  compact set.
  \item In each sector, $S_j,$ the operator $L$ agrees with an operator
    $L_j$ for which we can numerically construct the kernel of the limiting
    absorption resolvent $(L_j+k^2+i0)^{-1}.$ 
 \item  In the subsets $\{S_j:\: j=1,\dots,N\}$ we express the solution as a sum of an incoming solution
   and an outgoing solution $u_j=v^{\out}_j+v_j^{\In}$ with
   $(L_j+k^2)v^{\out}_j=(L_j+k^2)v^{\In}_j=0.$  The solutions $\{v_j^{\In}\}$
   are data encoding the incoming field. Each outgoing component
   $v^{\out}_j$ is expressed as a sum of layer potentials over the $\pa S_j$ using the
   kernel of $(L_j+k^2+i0)^{-1}$  and its normal derivative
   \item The global solution is found by
   imposing transmission boundary conditions to ensure that the function defined by
   \begin{equation}
     u^{\tot}({\bx})=u_j({\bx})\text{ for }{\bx}\in S_j, j=0,1,\dots,N,
   \end{equation}
   belongs to $H^2_{\loc}(\bbR^2).$ This is the requirement the $u_{\tot}({\bx}),\nabla
   u_{\tot}({\bx})$ are continuous, in the trace sense, across the boundaries of the sets
   $\{\pa S_j:\: j=0,\dots,N\}.$
  \end{enumerate}
  \item Using the integral representations, $G_j({\bx};{\by}),$ for $(L_j+k^2+i0)^{-1},$  the
    transmission boundary conditions are re-expressed as systems of integral
    equations on $\cup_{j}\pa S_j.$
    \item These integral equations are solved using the ``coordinate complexification'' method.
\end{enumerate}

 We call scatterers that satisfy 1(b) \emph{decomposable}.
The relationship between the data, $\{v_j^{\In}\}$  and
   the data $f$ in~\eqref{eqn1.2} depends on the problem one is solving. For
   example, if the
   incoming field comes from a point source, $G_j({\bx};{\by}_0)$  in $S_j,$ then $f$
   would usually be defined by choosing a cut-off $\psi$ equal 1 near to ${\by}_0$
   with support in $S_j$, and setting 
   \begin{equation}
     f=(L_j+k^2)(\psi({\bx})G_j({\bx};{\by}_0)).
   \end{equation}
The incoming data for the transmission problem is obtained by restricting
$G_j({\bx};{\by}_0)$ and $\pa_{\bn_{\bx}}G_j({\bx};{\by}_0)$ to $\pa S_j.$ Similar considerations
apply to other kinds of data, see~\cite{EpWG2023_1}.
 In a series of papers,~\cite{EpWG2023_1,EpWG2023_2,EpMaSR2023}, it
is \emph{proved} that this method can be used to solve the scattering problem if the
potential has the form
\begin{equation}\label{eqn1.6}
  q(x_1,x_2)=q_l(x_2)\chi_{(-\infty,0]}(x_1)+q_r(x_2)\chi_{[0,\infty)}(x_1).
\end{equation}
In this paper we implement the method described above for this example, and show
that it is also successful in much greater generality.

Key to the success of this program is the possibility of constructing kernels
for the resolvents $(L_j+k^2+i0)^{-1}.$  In general, constructing these kernels may be a very
difficult problem. However, if the potential has enough symmetry, then this construction often
becomes a practically solvable problem. 
For example, let ${\bx}=({\bx}',{\bx}'')$ be coordinates for
$\bbR^d$ and assume $L=\Delta+q({\bx}'),$ where $q$ is a compactly supported
function of the ${\bx}'$-variables. By taking Fourier transform in the remaining
variables we obtain the operators
\begin{equation}
  \tL[\xi'']=\Delta_{{\bx}'}-|\xi''|^2+q({\bx}').
\end{equation}
If $q({\bx}')$ has compact support, then the explicit construction of
$(\tL[\xi'']+k^2+i0)^{-1}$ can often be done quite efficiently. If $|\xi''|>k,$
then $(\tL[\xi'']+k^2)^{-1}$ is invertible, except possibly for $|\xi''|$ in a
finite set. The L.A.P. resolvent $(L+k^2+i0)^{-1}$ can then be constructed using
Fourier synthesis. As the portions of the solutions corresponding to $|\xi''|>k$
are largely evanescent, only frequencies in a compact set are required for
numerical solutions with a specified accuracy.  This approach to the construction of the
L.A.P. resolvent has been used by many authors: for example it used for the case
of layered media in~\cite{santosa2001wave,Christiansen1998}, and for the unreduced 2-body
problem in~\cite{Vasy97}. Even having a discrete group of translational
symmetries can make the construction of $(L+k^2+i0)^{-1}$ a computationally
tractable problem, see~\cite{AgocsBarnett_2024}.

In \cite{LuLuQian} a method for solving integral equations on infinite lines is introduced that we call the \emph{complexification} method.  This approach involves deforming the contour of integration into the complex
plane, taking advantage of the fact that the various fundamental solutions have analytic
continuations to certain domains, and for some choices of incoming and outgoing
variables, the analytically continued kernel functions are exponentially decaying.  In~\cite{BonChaFli_2022,BonChaFliComp_2022,EGrHJR_2024, tambova2018adiabatic} similar methods are used to solve integral equations.  This idea is  quite similar to that used in the complex scaling method, as described in Chapter 4.5 of~\cite{DZ}.

In order for the complexification method to work, it is necessary that the data also has analytic
extension properties to appropriate domains,  with exponential decay as  functions of $y_1+iy_2.$ This is often the case for physically meaningful data. For
example if we assume the integral equation is evaluated along the $x_1$-axis and
take as incoming data the fundamental solution itself, with source off of this
line, then, as $y_1\to \pm \infty,$ the function extends analytically to the
sets $Q_+=\{(y_1+it):\:0<t,y_1\}$ and $Q_-=\{(y_1+it):\:0>y_1,t\}.$ If $a+ib$
belongs to one of these quadrants, then
\begin{equation}
  H^{(1)}_0(k\sqrt{x_0^2+(t(a+ib)-y_0)^2})\sim \frac{e^{ikt(|a|+i|b|)}}{\sqrt{t(|a|+i|b|)}}[a_0+O(t^{-1})],
\end{equation}
decays exponentially. The square root is assumed to be defined on
$\bbC\setminus (-\infty,0),$ with positive values on $(0,\infty).$ As explained
in~\cite{EpWG2023_2}, this  implies that the solution to the integral equation has an
asymptotic expansion, and the same analytic extension and exponential decay
properties as the data and the kernel function. A simple example, with a
rigourous analysis of the complexification method, is given in~\cite{EGrHJR_2024}.

In Section~\ref{sec:single_guide} we use complexification to construct the L.A.P. resolvent kernel, $G_q(\bx;\by),$ in the case of a piece-wise constant, bi-infinite 2-dimensional wave-guide. We suppose that the wave-guide lies in the region $-d<x_2<d,$ so that
\begin{equation}
    q(x_2) = (k_1^2-k^2)\chi_{\{|x_2| <d\}}(x_2).
\end{equation}
In this case, we introduce the local wave-number
\begin{equation}
    k(\bx)=\begin{cases}
        k&\text{ if }|x_2|>d,\\
        k_1&\text{ if }|x_2|<d.
    \end{cases}
\end{equation}
In each region we express the kernel as a perturbation of the appropriate free space fundamental solution, $G_q(\bx;\by)=G_{k(\bx)}(\bx;\by)+G^{\out}_q(\bx;\by).$ The perturbation is found by using a boundary integral equation along $\{x_2=\pm d\}$ to solve for $G^{\out}_q,$  which is determined by transmission boundary conditions across the lines $\{x_2=\pm d\},$ and the outgoing radiation condition. This method is explained in detail in Section~\ref{sec:single_guide}.  With these Green's functions in hand, one recasts the scattering problem as a transmission problem.

Most of this paper is devoted to a detailed description of this approach applied
to potentials of the form given in~\eqref{eqn1.6}. We explain the construction
of the kernels of the limiting absorption resolvents
$(\Delta+q_{l,r}(x_1)+k^2+i0)^{-1}.$ A practical construction was given
in~\cite{EpWG2023_1}; here we use a different method.  We then set up the integral
equations along the line $\{x_1=0\},$ and explain how to use complexification to
solve this system of equations. This is followed by many examples, which are
solved using variants  of our method. For scattering problems of this sort, one expects the
incoming energy to equal the outgoing energy. We give a proof that this is the
case for open wave guides in an Appendix. Computing these quantities, we show
that our solutions have this property to many digits of accuracy.

While the method described above has only been rigorously justified for potentials of the
form in~\eqref{eqn1.6}, we show, in a series of examples, that the method
actually provides accurate answers in more
interesting and complicated geometries. These include, among others, semi-infinite wave guides
with compact perturbations at the end, and a pair of open wave-guides meeting at
various angles. These examples demonstrate the wider applicability of this approach.
\bigskip

\centerline{\bf Acknowledgements}
The authors would like to express our deep gratitude to Shidong Jiang for introducing us to the ideas underlying the complexification method, and encouraging us to use this approach to solve wave guide problems. We would also like to thank Jeremy Hoskins and Manas Rachh for their  very important contributions to many aspects of this project. Finally, we would  like to thank the American Institute of Mathematics and, in particular, John Fry for hosting us on Bock Cay during the SQuaREs program in May 2024. That meeting helped to launch this project in earnest.

\section{Review of Matched Wave-guides}
In this paper we first focus on the case of 2 semi-infinite wave-guides meeting along
a common perpendicular line, which is the decomposable problem treated
in~\cite{EpWG2023_1}.  This is modeled using a potential like that given
in~\eqref{eqn1.6}. For simplicity we consider the case of piecewise constant potentials,
\begin{equation}
  q_{l,r}(x_2)=(k_{1;l,r}^2-k^2)\chi_{[-d_{l,r},d_{l,r}]}(x_2),
\end{equation}
 though the more general case of bounded piecewise, continuous potentials can
 also be handled using this method. We  usually assume that the wave numbers, $k_{1;l,r},$ inside
 the channel are larger than that in the surrounding medium. This is needed 
 for ``wave-guide modes'' to exist, and is the case of principal interest in
 opto-electronics.

\subsection{Wave-guide modes}\label{ss.WGmodes}
For an operator of this type there are finitely many solutions $v_{l,r;n}(x_2)\in L^2(\bbR)$ to
\begin{equation}\label{eqn2.2.002}
    (\pa_{x_2}^2 + q_{l,r}(x_2)) v_{l,r;n}(x_2) = \lambda_{l,r;n}^2v_{l,r;n}(x_2),
    n=1,\dots,N_{l,r},
\end{equation}
where $0<\lambda_{l,r;n}^2<k_{1;l,r}^2-k^2.$ 
In fact, such a solution is easily seen to be $O(e^{-\lambda_{l,r;n}|x_2|})$ as
$|x_2|\to\infty.$ The frequencies for which such a solution exist can be found
by solving a simple transcendental equation. If the channel has width $2d$ with
 interior wave number  $k_{1;l,r},$ and  exterior wave number  $k,$ then
the wave guide energies are found by solving
\begin{equation}\label{eq:mode_eqn}
  \tan 2d\sqrt{k_{1;l,r}^2-k^2-\lambda^2}=\frac{2\lambda\sqrt{k_{1;l,r}^2-k^2-\lambda^2}}{k_{1;l,r}^2-k^2-2\lambda^2}.
\end{equation}

If we let $\xi_{l,r;n}=\sqrt{k^2+\lambda_{l,r;n}^2}>0,$ then we
have that
\begin{equation}
  (\Delta+q_{l,r}+k^2)e^{\pm i\xi_{l,r;n}x_1}v_{l,r;n}(x_2)=0;
\end{equation}
these are the \emph{wave-guide modes}.  They are not $L^2$-eigenfunctions in 2
dimensions as they do not decay in the $x_1$-direction, but they are strongly
localized within the $\supp q_{l,r}.$ The functions
\begin{equation}
  \{e^{i\xi_{l,r;n}x_1}v_{l,r;n}(x_2):n=1,\dots,N_{l,r}\}
\end{equation}
are the rightward moving wave-guide modes and
\begin{equation}
  \{e^{-i\xi_{l,r;n}x_1}v_{l,r;n}(x_2):n=1,\dots,N_{l,r}\}
\end{equation}
are the leftward moving wave-guide modes.

As noted in our earlier work, two semi-infinite wave-guides meeting along a
common perpendicular line is the simplest case for which it is not possible to
explicitly compute the fundamental solution, and we instead follow the outline
above and rephrase the scattering problem as a transmission problem along the
common perpendicular line $\{x_1=0\}.$ Typical incoming data for such a
wave-guide structure is given by a sum of incoming wave-guide modes, i.e.
\begin{equation}
  v^{\In}_l=\sum_{n=1}^{N_l}c_ne^{i\xi_{l,r;n}x_1}v_{l,r;n}(x_2),\text{ for }x_1<0,
\end{equation}
though other types of incoming data are allowed, see Section 6 of~\cite{EpWG2023_1}.
The scattered field is a sum of outgoing wave-guide modes as well as outgoing
radiation. 
\subsection{IE for matched Wave-guides}
The first step in replacing the scattering problem with a transmission problem
is the construction of outgoing fundamental solutions for the operators
$(\Delta+q_{l,r}+k^2),$ whose kernels we denote by $G_{q_{l,r}}(\bx;\by).$ We
assume that there are incoming fields $v_{l,r}^{\In}$ defined in $\pm x_1>0,$ which
solve the equations
\begin{equation}
  (\Delta+q_{l,r}+k^2)v_{l,r}^{\In}=0\text{ where }\pm x_1>0.
\end{equation}
The total field is given by
\begin{equation}
  u^{\tot}(x_1,x_2)=
  \begin{cases}
    v_{l}^{\In}(x_1,x_2)+v_l^{\out}(x_1,x_2)&\text{ where }x_1<0,\\
         v_{r}^{\In}(x_1,x_2)+v_r^{\out}(x_1,x_2)&\text{ where }x_1>0;
  \end{cases}
\end{equation}
the $v_{l,r}^{\out}$ are ``outgoing'' solutions to $
(\Delta+q_{l,r}+k^2)v_{l,r}^{\out}=0$ where $\pm x_1>0.$ These solutions are
determined by the transmission boundary conditions
\begin{equation}\label{eqn2.9.001}
  \begin{split}
     v_{l}^{\In}(0^-,x_2)+v_l^{\out}(0^-,x_2)&=
     v_r^{\In}(0^+,x_2)+v_r^{\out}(0^+,x_2)\\
       \pa_{x_1}v_{l}^{\In}(0^-,x_2)+\pa_{x_1}v_l^{\out}(0^-,x_2)&= \pa_{x_1}v_r^{\In}(0^+,x_2)+\pa_{x_1}v_r^{\out}(0^+,x_2).
  \end{split}
\end{equation}

Provisionally, the outgoing condition is imposed by representing $v^{\out}_{l,r}$
in terms of layer potentials constructed out of the outgoing fundamental solutions, for $\pm x_1>0:$
\begin{equation}\label{eq:u_rep}
\begin{split}
  v^{\out}_{l,r}(\bx)&=\cS_{q_{l,r}}[\tau](\bx) + \cD_{q_{l,r}}[\sigma](\bx) \\
  &=\int_{-\infty}^{\infty}G_{q_{l,r}}(\bx;0,y_2)\tau(y_2)dy_2+
  \int_{-\infty}^{\infty}\pa_{y_1}G_{q_{l,r}}(\bx;0,y_2)\sigma(y_2)dy_2.
  \end{split}
  \end{equation}
  
As shown in~\cite{EpWG2023_1,EpWG2023_2,EpMaSR2023}, for appropriate incoming
data, the functions $v^{\out}_{l,r}$ define solutions that satisfy the physically
meaningful outgoing radiation conditions. The boundary conditions
in~\eqref{eqn2.9.001} lead to integral equations along the line $\{x_1=0\},$
\
\begin{equation}
  \begin{split}
   \lim_{x_1\to
      0^-}&\left[\int_{-\infty}^{\infty}(G_{q_{l}}(\bx;0,y_2)\tau(y_2)+\pa_{y_1}G_{q_{l}}(\bx;0,y_2)\sigma(y_2))dy_2\right]-\\
   & \lim_{x_1\to
     0^+}\left[\int_{-\infty}^{\infty}(G_{q_{r}}(\bx;0,y_2)\tau(y_2)+\pa_{y_1}G_{q_{r}}(\bx;0,y_2)\sigma(y_2))dy_2\right]=\\
   &v_r^{\In}(0,x_2)-v_l^{\In}(0,x_2);
     \end{split}
\end{equation}
    \begin{equation}
  \begin{split}
    \lim_{x_1\to
      0^-}\pa_{x_1}&\left[\int_{-\infty}^{\infty}(G_{q_{l}}(\bx;0,y_2)\tau(y_2)+\pa_{y_1}G_{q_{l}}(\bx;0,y_2)\sigma(y_2))dy_2\right]-\\
   & \lim_{x_1\to
      0^+}\pa_{x_1}\left[\int_{-\infty}^{\infty}(G_{q_{r}}(\bx;0,y_2)\tau(y_2)+\pa_{y_1}G_{q_{r}}(\bx;0,y_2)\sigma(y_2))dy_2\right]=\\
    &\pa_{x_1}v_r^{\In}(0,x_2)-\pa_{x_1}v_l^{\In}(0,x_2)
  \end{split}
    \end{equation}

    Using standard jump conditions for layer potentials, and a detailed analysis
    of the kernels $G_{q_{l,r}}(\bx;\by)$ it is shown in~\cite{EpWG2023_1} that
    these equations take the rather simple form
    \begin{equation}\label{eqn2.13.001}
      \left(
      \begin{matrix}
        \Id&D\\C&\Id
      \end{matrix}\right) \left(
      \begin{matrix}
        \sigma\\\tau
      \end{matrix}\right)=\left(
      \begin{matrix}
        v_r^{\In}(0,x_2)-v_l^{\In}(0,x_2)\\\pa_{x_1}[v_r^{\In}(0,x_2)-v_l^{\In}(0,x_2)]
      \end{matrix}\right).
    \end{equation}
    The operator on the left hand side is a Fredholm operator of index zero acting on the spaces
    $\cC_{\alpha}(\bbR)\oplus\cC_{\alpha+\frac 12}(\bbR)$ for $0<\alpha<\frac
    12$ of weighted continuous functions. The norm  on $\cC_{\beta}(\bbR)$ is defined by
    \begin{equation}
      |f|_{\beta}=\sup\{(1+|x|)^{\beta}|f(x)|:\:x\in\bbR\}.
    \end{equation}
    In~\cite{EpMaSR2023} it is shown that the null spaces of these equations are
    trivial. This shows that these equations are uniquely solvable for arbitrary data
    $(g,h)\in\cC_{\alpha}(\bbR)\oplus\cC_{\alpha+\frac 12}(\bbR).$  In order for the corresponding solution to the transmission problem to satisfy the physically meaningful radiation conditions given in~\cite{EpMaSR2023}, the data must also satisfy appropriate radiation conditions.

    In the following sections we give numerical approximations for the kernels
    of $C$ and $D;$ in~\cite{CmplxWG2} we show that they continue analytically to certain subsets
    of $\bbC\times\bbC.$ For data $(g,h)$ that also has analytic continuations
    we can deform the contour on which we are solving the integral equation.
    The analytic continuations of the data and the solutions
    $(\sigma,\tau)$ are exponentially decaying.  This allows us to solve the
    integral equations for $\{(\sigma(x_2),\tau(x_2)):\: x_2\in [-L,L]\},$ for
    any fixed $L>0$ and with any specified accuracy, by replacing the infinite
    domain implicit in~\eqref{eqn2.13.001} with a finite interval on the
    deformed contour. See the end of Section~\ref{ss.coord_complx}.
    
    \subsection{Asymptotics of Solutions and the Scattering Matrix}\label{ss.scat_matrix.003}
    Of special interest in problems connected to open wave-guides is the
    scattering matrix, which relates the leading incoming part of the solution
    to the leading outgoing part. Suppose $u$ is a solution to
    $(\Delta+q+k^2)u=0,$ where $q$ is a potential of the type defined
    in~\eqref{eqn1}.  Let $S^1$ be identified with the boundary of the radial
    compactification of $\bbR^2,$ with $\{\omega_1,\dots,\omega_N\}$ the
    endpoints of the channels at infinity. Assume that $u$ has an asymptotic expansion at
    infinity of the form
    \begin{equation}\label{eqn2.15.001}
      u(r\omega)\sim
      \frac{e^{ikr}}{\sqrt{r}}a_+(\omega)+\frac{e^{-ikr}}{\sqrt{r}}a_{-}(\omega)+o(r^{-\frac 12}).
    \end{equation}
    The coefficients $a_{\pm}\in L^{\infty}(S^1)$ are smooth, functions on
    $S^1\setminus\{\omega_1,\dots,\omega_N\}.$

    Along the channels, $u$ satisfies:
    \begin{equation}\label{eqn2.16.001}
      u(\bx)=\sum_{n=1}^{M_j}[a_{n}^je^{i\xi_n^j\langle\bx,\omega_j\rangle}+b_{n}^je^{-i\xi_n^j\langle\bx,\omega_j\rangle}]
      v_n^j(\langle\bx,\omega^{\bot}_j\rangle)
      +O(|\langle\bx,\omega_j\rangle|^{-\frac 12}),
    \end{equation}
with $\langle\omega_j,\omega_j^{\bot}\rangle=0$ and
$\det(\omega_j\,\omega_j^{\bot})=1.$ As shown in~\cite{EpMaSR2023}, this is true
of solutions with a purely radial scattering wave-front set.  The scattering
matrix relating the incoming asymptotics to the outgoing asymptotics is defined
by
\begin{multline}
  S(k):\left(\sqrt{k}a_-(\omega),\left[\sqrt{\xi_n^j}b_n^j:\: n=1,\dots,M_j,\, j=1,\dots,
    N\right]\right)\mapsto\\
  \left(\sqrt{k}a_+(\omega),\left[\sqrt{\xi_n^j}a_n^j:\: l=1,\dots,M_j,\, j=1,\dots, N\right]\right).
\end{multline}

In Appendix~\ref{App1} we prove the following conservation law, which is equivalent to
the unitarity of the scattering matrix.
\begin{theorem}\label{scat_unit_thm}
  If $u$ is solution to $(\Delta+q+k^2)u=0,$ with asymptotic expansions like
  those in~\eqref{eqn2.15.001} and ~\eqref{eqn2.16.001}, then
  \begin{equation}\label{eqn2.18.001}
    k\int_{S^1}|a_-(\omega)|^2d\omega+\sum_{j=1}^N\sum_{n=1}^{M_j}\xi_n^j|b_n^j|^2=
     k\int_{S^1}|a_+(\omega)|^2d\omega+\sum_{j=1}^N\sum_{n=1}^{M_j}\xi_n^j|a_n^j|^2.
  \end{equation}
\end{theorem}

In our applications, the only incoming part of the solution comes from the
data, and we construct an outgoing solution $u^{\out}.$ In the proof of 
Theorem~\ref{scat_unit_thm} we show that
\begin{equation}\label{eqn2.19.001}
  \lim_{R\to\infty}\Im\left[\int_{0}^{2\pi}\overline{u^{\out}(R\omega)}\pa_ru^{\out}(R\omega)Rd\omega
    \right]= k\int_{S^1}|a_+(\omega)|^2d\omega+\sum_{n=1}^{M_j}\xi_n^j|a_n^j|^2.
\end{equation}
In~\cite{EpMaSR2023} we also show that the coefficients $\{|a_n^j|\}$ can also be found by
computing the limits
\begin{equation}
|a_n^j|=  \lim_{x_1\to\infty}\left|\int_{-\infty}^{\infty}u^{\out}(x_1,x_2)v_n^j(x_2)dx_2\right|.
\end{equation}
In~\cite{EpWG2023_1} we show that for the special case of 2 semi-infinite wave guides
meeting along a common perpendicular line
\begin{equation}\label{eqn2.22.001}
  \begin{split}
   & a_n^re^{i\xi_n^rx_1}=\int_{-\infty}^{\infty}v_{r}^{\out}(x_1,x_2)v_n^r(x_2)dx_2\text{
      for any } x_1>0,\\
    & a_n^le^{-i\xi_n^lx_1}=\int_{-\infty}^{\infty}v_{l}^{\out}(x_1,x_2)v_n^l(x_2)dx_2\text{
      for any } x_1<0.
  \end{split}
\end{equation}

The solutions we construct to the scattering problem certainly satisfy the PDE
everywhere, and belong to $H^2_{\loc}(\bbR^2).$ The only real question is
whether the solutions we construct are truly outgoing.  As there is essentially no
data for which there is a gold standard for numerical solutions of open wave
guide scattering problems, these observations give a useful way to assess the
quality of the solutions found by a numerical algorithm. The left hand side
of~\eqref{eqn2.18.001} is computable from the incoming field, which is specified
as data. Equations~\eqref{eqn2.19.001}--\eqref{eqn2.22.001} provide a
computationally effective means to approximate the right hand side
of~\eqref{eqn2.18.001}. While a very small difference between the two sides does
not prove the pointwise accuracy of a solution, it does provide strong evidence
that the solution we obtain satisfies the outgoing radiation condition.

From~\eqref{eqn2.15.001} we also see that the computed solution $u^{\out}$ has an expansion
\begin{equation}\label{eqn2.24.004}
   \frac{R}{k}\overline{u^{\out}(R\omega)} \pa_r u^{\out}(R\omega)
   =i[|a_+(\omega)|^2-|a_-(\omega)|^2]-2\Im[e^{2ikR}a_+(\omega)\overline{a_-(\omega)}]+O(R^{-1}),
\end{equation}
provided $\omega\notin\{\omega_1,\dots,\omega_N\}.$
Evaluating this for a large $R$ gives a pointwise estimate on the failure of the solution to be outgoing. In particular, we should find that 
$$\Re\overline{u^{\out}(R\omega)} \pa_r u^{\out}(R\omega)=O(R^{-2}).$$

\section{Green's function for a single wave-guide}\label{sec:single_guide}

In order to solve~\eqref{eqn2.13.001} we must be able to compute a bi-infinite wave-guide Green's function. There is  a vast literature on this subject (see, for example~\cite{bruno2017windowed,bruno2016windowed,cai2013computational,okhmatovski2004evaluation,wang2019fast,zhang2020exponential,paulus2000accurate,taraldsen2005complex}). One popular method is to compute the inverse Fourier transform of the L.A.P. resolvent kernels of $\pa_{x_2}-\xi^2+k^2+q(x_2)$. The resulting formulas are also referred to as Sommerfeld integral formulas and are used in \cite{cai2013computational,okhmatovski2004evaluation,wang2019fast,zhang2020exponential,paulus2000accurate}. An alternative approach is the windowed Green's function approach developed in~\cite{bruno2017windowed,bruno2016windowed}. In this section, we briefly present a new method based on the coordinate complexification technique developed in~\cite{EGrHJR_2024}. Indeed, any method may be used to find the Green's function without changing the subsequent steps needed to solve the scattering problem. We choose this method so that we can use the same discretization techniques in~\eqref{eq:wave_guide_IE},~\eqref{eqn2.13.001} and the other integral equations in Section~\ref{sec:other_ex}. 

In this paper we consider piece-wise constant wave-guides, which corresponds to the choice of~$q(x_2) =  (k_1^2-k^2) \chi_{[-d,d]}(x_2)$. As we mentioned in the introduction, we write a bi-infinite wave-guide Green's function as $G_q(\bx;\by)=G_{k(\bx)}(\bx;\by)+G^{\out}_q(\bx;\by).$ It is not hard to see that~$G_q^{\out}$ is the ``outgoing" solution of the transmission problem
\begin{equation}
\begin{cases}
        (\Delta + k^2+q(x_2)) G_q^{\out}(\bx;\by)= 0\quad \text{in}\quad \bbR^2\setminus\{|x_2|=d\},\\
[[G_q^{\out}]]_{|x_2|=d}  = -[[G_{k(\bx)}]]_{|x_2|=d} ,\\  [[\partial_{x_2} G_q^{\out}]]_{|x_2|=d}  = -[[\partial_{x_2} G_{k(\bx)}]]_{|x_2|=d} .
\end{cases}\label{eq:G_2_eq}
\end{equation}
To numerically solve for~$G_q$, it is therefore enough to be able solve transmission problems of this form. In the remainder of this section we discuss a numerical solver for transmission problems of this type.  The decomposition used here differs from that in~\cite{EpWG2023_1,EpWG2023_2} where we write $G_{q}(\bx;\by)=G_k(\bx;\by)+\tilde{G}^{\out}_q(\bx;\by).$

Before we do, however, we observe that this decomposition can easily be used to compute derivatives of~$G_q$: the derivatives of~$G_{k(\bx)}$ can easily be computed and the derivatives of~$G_q^{\out}$ can be computed by differentiating the solution or the boundary data of~\eqref{eq:G_2_eq}.

We finish this preliminary discussion by showing how the Sommerfeld integral representation leads to a solution, in principle, to this problem that clarifies how such a representation ``finds" the correct wave-guide modes.  In our actual computations we do not use the Fourier representation, but it becomes clearer why this approach works. We begin with the fact that the outgoing fundamental solution is obtained as the $\lim_{\delta\to 0^+}(\Delta+k^2(x_2)+i\delta)^{-1},$ and the Sommerfeld integral representation for the free space fundamental solution
\begin{equation}
    G_{\sqrt{k^2+i\delta}}(\bx;\by)=
    \frac{i}{4\pi}\int_{-\infty}^{\infty}\frac{
    e^{i|x_2-y_2|\sqrt{k^2+i\delta-\xi^2}}e^{i(x_1-y_1)\xi}d\xi}{\sqrt{k^2+i\delta-\xi^2}}.
\end{equation}

Assume that the pole of~$G_{\sqrt{k^2+i\delta}}$ is $\by=(0,y_2),$ where $y_2>d;$ similar considerations apply in all cases. We write the kernel, $R_{\delta}(\bx;\by),$ for $(\Delta+k(x_2)^2+i\delta)^{-1}$ as 
$$R_{\delta}(\bx;\by)=G_{\sqrt{k^2+i\delta}}(\bx;\by)\cdot\chi_{[d,\infty)}(x_2)+e_{\delta}(\bx;\by).$$ 
The correction term, $e_{\delta}$ satisfies $(\Delta_{\bx}+k(x_2)^2+i\delta)e_\delta=0$ in the complement of the lines $\{x_2=\pm d\},$ and is determined by transmission boundary conditions across these lines.

In each of the three regions we can express $e_{\delta}$ as sums of single and double layer potentials:
\begin{equation}
    \begin{split}
        &e_{\delta}(\bx;\by)=\\&\int_{-\infty}^{\infty}
        \left[G_{\sqrt{{k}^2+i\delta}}(\bx;t,\pm d )\tau_{\delta}^{\pm}(t)
       \mp \pa_{y_2}G_{\sqrt{{k}^2+i\delta}}(\bx;t,\pm d)\sigma_{\delta}^{\pm}(t)\right]dt,\text{ for }\pm x_2>d,\\
       &\phantom{ee}\\
        &e_{\delta}(\bx;\by)=\\&\int_{-\infty}^{\infty}
        \left[G_{\sqrt{k_1^2+i\delta}}(\bx;t,d)\tau_{\delta}^{+}(t)
       + \pa_{y_2}G_{\sqrt{k_1^2+i\delta}}(\bx;t,d)\sigma_{\delta}^{+}(t)\right]dt+\\
       &\int_{-\infty}^{\infty}
        \left[G_{\sqrt{k_1^2+i\delta}}(\bx;t,-d)\tau_{\delta}^{-}(t)
       - \pa_{y_2}G_{\sqrt{k_1^2+i\delta}}(\bx;t,-d)\sigma_{\delta}^{-}(t)\right]dt,\text{ for }|x_2|<d.
    \end{split}
\end{equation}
For $\delta>0$ it is easy to see that the densities $\{\sigma_{\delta}^{\pm},\tau_{\delta}^{\pm}\}$ are decaying and have well behaved Fourier transforms. The densities also depend on $\by,$ but to simplify the notation, we omit this argument. 

Using the Sommerfeld formula we can re-express $e_{\delta}(\bx;\by)$ in terms of the Fourier transforms of the densities, $\{\hsigma_{\delta}^{\pm},\htau_{\delta}^{\pm}\}.$ For example, if $x_2>0,$ then we have
\begin{equation}\label{eqn1.19.004}
    e_{\delta}(\bx;\by)=
   \frac{i}{4\pi}\int_{-\infty}^{\infty}\frac{
    e^{i(x_2-d)\sqrt{{k}^2+i\delta-\xi^2}}e^{ix_1\xi}
    \left[\htau^+_{\delta}(\xi)+i\sqrt{{k}^2+i\delta-\xi^2}\hsigma^+_{\delta}(\xi)\right]d\xi}{\sqrt{{k}^2+i\delta-\xi^2}}. 
\end{equation}
Using these representations for $e_{\delta}(\bx;\by)$ and $\pa_{x_2}e_{\delta}(\bx;\by)$ where $x_2=\pm d,$ with $f(\pm d^{\pm})=\lim_{x\to\pm d^{\pm}}f(x),$ the transmission boundary conditions are:
\begin{equation}
\begin{split}
    &e_{\delta}(x_1,d^+;\by)+G_{k}(x_1,d;\by)=e_{\delta}(x_1,d^-;\by),\\ 
    &e_{\delta}(x_1,-d^+;\by)=e_{\delta}(x_1,-d^-;\by),\\
      &\phantom{ee}\\
    &\pa_{x_2}e_{\delta}(x_1,d^+;\by)+\pa_{x_2}G_{k}(x_1,d;\by)=\pa_{x_2}e_{\delta}(x_1,d^-;\by),\\ &\pa_{x_2}e_{\delta}(x_1,-d^+;\by)=\pa_{x_2}e_{\delta}(x_1,-d^-;\by).
    \end{split}
\end{equation}

In the Fourier representation these take the form of  a family of $4\times 4$ linear systems:
\begin{equation}
    A_{\delta}(\xi)\fv_{\delta}(\xi)=\fd_{\delta}(\xi),\text{ for }\xi\in\bbR.
\end{equation}
Here $\fv_{\delta}(\xi)=(\hsigma_{\delta}^{+}(\xi),\htau_{\delta}^{+}(\xi),\hsigma_{\delta}^{-}(\xi),\htau_{\delta}^{-}(\xi))^t,$ and $\fd_{\delta}(\xi)$ is obtained from the Sommerfeld representation of the data $G_{\sqrt{{k_1}^2+i\delta}}(\bx;\by).$  It is well known that these matrices are invertible for $\delta>0.$ Let $A_{0}(\xi)=\lim_{\delta\to 0^+}A_{\delta}(\xi);$ this matrix is invertible except for frequencies $\{\xi_n\}$ for which
\begin{equation}
    (\pa_{x_2}^2+k^2(x_2)-\xi_n^2)v_n(x_2)=0
\end{equation}
has an $L^2$-solution. There are finitely many such frequencies, which lie in $[{k},k_1);$ these correspond to \emph{wave-guide modes} for the bi-infinite wave-guide, which are discussed in Section~\ref{ss.WGmodes}.

In fact the $\det A_{0}(\xi)$ has simple poles at $\{\pm\xi_n\},$ which complicates the evaluation of $\lim_{\delta\to 0^+}\fv_{\delta}(\xi)=\lim_{\delta\to 0^+} A_{\delta}^{-1}(\xi)\fd_{\delta}(\xi).$ The simplest thing to do is to assume that $A_0(k)\neq 0$ and deform the contour of integration in~\eqref{eqn1.19.004}, e.g.  to a contour that replaces intervals around $\{\pm \xi_n\}$ with semicircles centered on these points in either the upper or lower half planes, and then let $\delta\to 0^+.$  Using the semi-circles in the upper half plane produces right-ward moving wave-guide modes, and those in the lower half plane, left-ward moving wave-guide modes. This is essentially the construction used in~\cite{EpWG2023_1}.  Since the poles of $\det A_{0}(\xi)$ are simple, the limiting singularity of $A_0^{-1}(\xi)$ takes the form 
$$\lim_{\delta\to 0^+}\frac{a_n}{\xi-\xi_n+i\delta}=a_n\left[\PV\frac{1}{\xi-\xi_n}-i\pi\delta(\xi-\xi_n)\right],$$
as tempered distributions, which gives an alternate way to analyze the limiting densities, which are themselves tempered distributions.

\subsection{Integral equation formulation}
To solve~\eqref{eq:G_2_eq}, we first write it in a more generic form. We let~$\gamma$ be the union of the lines~$x_2=\pm d$, which are the discontinuities of~$q$. We shall also let~$\bs n$ be the unit normal to~$\gamma$, defined to point away from the region with wavenumber~$k_1$. 
We will look for the solution~$v$ of
\begin{equation}
\begin{cases}
        (\Delta +k^2+ q(x_2)) v= 0\quad \text{in}\quad \bbR^2\setminus\gamma\\
[[v]]_{\gamma}  = -[[v^{\In}]]_{\gamma}, \quad [[\partial_{\bs n} v]]_{\gamma}  = -[[\partial_{\bs n}  v^{\In}]]_{\gamma}
\end{cases},\label{eq:Transmission_PDE}
\end{equation}
plus appropriate radiation conditions. Clearly~\eqref{eq:G_2_eq} is of this form with a particular choice of~$v^{\In}$.

In order to numerically solve~\eqref{eq:Transmission_PDE} we convert it into an integral equation supported on~$\gamma$. To this end, we introduce the single and double layer potentials
\begin{equation}
    \cS_{k}[\rho](\bx) = \int_{\gamma} G_k(\bx;\by)\rho(\by) d\by \quad \text{and}\quad \cD_{k}[\mu](\bx) = \int_{\gamma} \partial_{\bs n_{\by}}G_k(\bx;\by)\mu(\by) d\by.
\end{equation}
The operators~$\cS_{k_1}$ and~$\cD_{k_1}$ are defined similarly. As these layer potentials are integrals are over unbounded domains and the kernels decay quite slowly, a rigorous study of these operators would require careful consideration to ensure that the integrals exist. In this paper however, we proceed formally, and assume that the integrals make sense.

Under this assumption, the function
\begin{equation}
    v(\bx) = \begin{cases}
        \cD_{k}[\mu](\bx)+\cS_{k}[\rho](\bx) & \text{ if }|x_2|>d\\
        \cD_{k_1}[\mu](\bx)+\cS_{k_1}[\rho](\bx) & \text{ if }|x_2|<d
    \end{cases}\label{eq:wave_guide_rep}
\end{equation}
automatically satisfies the PDE in~\eqref{eq:Transmission_PDE}. We are therefore free to choose~$\rho$ and~$\mu$ so that~$v$ satisfies the transmission condition.

To enforce the transmission conditions, we introduce the normal derivatives of the layer potentials:
\begin{equation}
    \cS_{k}'[\rho](\bx) = \int_{\gamma} \partial_{\bs n_{\bx}}G_k(\bx;\by)\rho(\by) d\by \quad \text{and}\quad \cD'_{k}[\mu](\bx) = \int_{\gamma} \partial_{\bs n_{\bx}}\partial_{\bs n_{\by}}G_k(\bx;\by)\mu(\by) d\by,
\end{equation}
where the~$\cD'_k$ operator is interpreted in the Hadamard finite parts sense (see~\cite{kress}).

If we use the well known jump relations for the single and double layer potentials (see~\cite{kress}), then we find that~$v$ given by~\eqref{eq:wave_guide_rep}  satisfies the continuity conditions in~\eqref{eq:Transmission_PDE} if
\begin{equation}
    \lp \Id + \begin{pmatrix}
        \cD_{k} - \cD_{k_1} & \cS_{k} - \cS_{k_1}\\
        -(\cD'_{k} - \cD'_{k_1}) & -(\cS'_{k} - \cS'_{k_1})
    \end{pmatrix}\rp\begin{pmatrix}
        \mu \\ \rho
    \end{pmatrix} = \begin{pmatrix}
        -[[ v^{\In}]]_{\gamma} \\ [[\partial_{\bs n}  v^{\In}]]_{\gamma}
    \end{pmatrix} \quad\text{on}\quad \gamma.\label{eq:wave_guide_IE}
\end{equation}
We thus have that if~$\rho$ and~$\mu$ satisfy~\eqref{eq:wave_guide_IE}, then~$v$ defined by~\eqref{eq:wave_guide_rep}  solves~\eqref{eq:Transmission_PDE}.

\label{ss.coord_complx}
Due to the presence of wave-guide modes, the
solutions~$\rho$ and~$\mu$ do not decay as~$x_1\to \pm
\infty$. Therefore the integral
equation~\eqref{eq:wave_guide_IE} cannot be numerically truncated at any value of~$x_1$. We therefore use the coordinate complexification method introduced in the introduction. This is quite similar to the approach taken in~\cite{BBD_SIAM_2024}.

In this method, we deform the boundary~$\gamma$ to the curve
\begin{equation}
    \tilde\gamma = \{(x_1+i \psi(x_1),\pm d)\; |\; x_1\in \bbR \},
\end{equation}
where~$\psi$ is a monotonically increasing function supported on~$[-L,L]^C$, which goes to~$\pm\infty$ as~$x_1\to\pm\infty$. The choice of~$L$  determines the strip of~$\bbR^2$ on which we can recover the solution~$u$. We  choose the normal to~$\tilde\gamma$ to be~$\bs n= (0,\operatorname{sign}(x_2))$, which agrees with the normal to~$\gamma$ where they overlap.

The integral operators on~$\tilde\gamma$ are the same as those on~$\gamma$, except that they are based on the analytic continuation of the Green's function:
\begin{equation}
    \tilde G_k(\bx;\by) := \frac{i}4H_0^{(1)}(k r(\bx;\by)),
\end{equation}
where
\begin{equation}
    r(\bx;\by) := \sqrt{(x_1-y_1)^2+(x_2-y_2)^2}.
\end{equation}
Choosing~$\psi$ to be non-decreasing ensures that
\begin{equation}
    \Im[(x_1-y_1+i(\psi(x_1)-\psi(y_1))^2+(x_2-y_2)^2]\geq 0,
\end{equation}
for all~$x_1,y_1,x_2,y_2\in\bbR$. Therefore $r(\cdot;\by)|_{\tilde \gamma}$ does not cross the branch cut of the square root.

The advantage of this contour deformation comes from the fact that our kernels and densities are asymptotic to linear combinations of functions of the form
\begin{equation}
    e^{i\eta \sqrt{x_1^2}} (x_1^2)^{\alpha/2}
\end{equation}
for various choices of~$\eta$ and~$\alpha$. When functions of this form are analytically continued to~$\tilde\gamma$, they are bounded by the exponentially decaying functions
\begin{equation}
   C e^{- \eta |\Im x_1|}|x_1|^\alpha.\label{eq:exp_deacy}
\end{equation}
This exponential decay allows us to truncate~$\tilde\gamma$ when~$|\Im x_1|$ is large enough.

\subsection{Discretization}\label{sec:wave_discretize}

We now discuss a numerical method for solving the integral equation~\eqref{eq:wave_guide_IE} on the contour~$\tilde \gamma$. For the remainder of this paper, we take the complexification function to defined by $\psi(x_1) = 20 \operatorname{erfc}( (x_1+L+30)/5)-\operatorname{erfc}( -(x_1-L-30)/5)$, which is smooth and  zero to machine precision on~$[-L,L]$.

Since the densities decay exponentially on the non-real portions of~$\tilde\gamma$, we can safely truncate without a great loss of accuracy. Specifically, since the slowest oscillation occurs at wavenumber~$k$, there exists a constant~$C>0$ such that
\begin{equation}
    |\rho|,|\mu| < Ce^{-k|\Im x_1|}.
\end{equation}
The estimate in~\eqref{eq:exp_deacy} implies that we may truncate the contour when $|\psi(x_1)| k > \log \epsilon $, where~$\epsilon$ is a desired accuracy. In our experiments $k$ is $\sim 1$ and we choose the truncation~$\epsilon$ to be~$10^{-17}$, so we can truncate at~$\Re x_1=\pm(L+38)$.

On this truncated contour, the integral equation can be discretized using standard techniques. Specifically, we solve it using the chunkIE package~\cite{chunkIE}, which solves  integral equations using a modified Nystr\"om discretization. The contour~$\tilde\gamma$ is split into panels~$\gamma_i$, each of which is represented by 16th order Gauss-Legendre nodes. Overall, the contour is represented by nodes~$\bx_1,\ldots, \bx_{\nchnk}$. The unknowns~$\rho$ and~$\mu$ are represented by their values at the~$\bx_i$'s, which we denote by~$\rho_i$ and~$\mu_i$.

To evaluate integral operators at a target~$\bx$, we do integrals one panel at a time. If the target is far from the panel~$\gamma_i$, then we use the Gauss-Legendre quadrature rule, which is well suited to the smooth quadratures. When~$\bx$ is near but not on~$\gamma_i$, we use adaptive integration to integrate the interpolant of~$\rho$ and~$\mu$. When~$\bx$ is in~$\gamma_i$, we use generalized Gaussian quadrature~\cite{bremer2010nonlinear}. With the above quadrature rules, we can build a discrete linear system that enforces~\eqref{eq:wave_guide_IE} at each node~$\bx_i$.

We solve this discretized equation using the recursive skeletonization algorithm implemented in the \emph{Fast Linear Algebra in MATLAB} (FLAM) software library~\cite{ho2020flam}. This algorithm builds a fast direct solver for~\eqref{eq:wave_guide_IE}, which is particularly appropriate for this problem because we must solve the equation with a different right hand side every time we wish to evaluate~$G_q$. The cost of building this fast direct solver is $O(n^2)$, and once it has been built the cost of each solve is $O(n)$. Once we have solved~\eqref{eq:wave_guide_IE}, the solution~$u$ can be recovered through the representation~\eqref{eq:wave_guide_rep} and the previously discussed quadrature rules.

\subsection{Numerical verification}\label{sec:test_single}
\begin{figure}
    \centering
    \includegraphics[width=.6\textwidth]{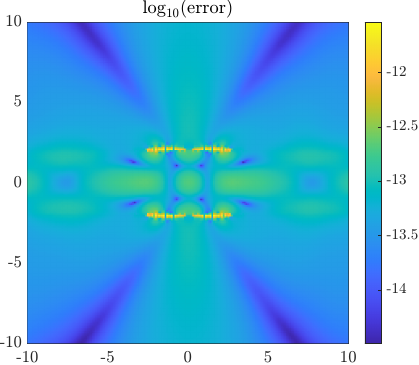}
    \caption{The analytic solution test for our single wave-guide solver.
    We see that the solver is accurate to at least the adaptive integration tolerance~$10^{-12}$ everywhere.}
    \label{fig:single-ana_test}
\end{figure}

We test our solver with~$k=1,k_1=2,$ and~$d=2$, and begin with an analytic solution test by picking the incoming field~$v^{\In}$ such that the true solution is known. Specifically, if we let~$\bx_0=(0,1)$, then it is not hard to see that 
\begin{equation}
    v^{\In}(\bx) = \begin{cases}
        \tilde G_k(\bx;\bx_0) & \text{ if }|x_2|>d\\
        0 & \text{ if }|x_2|<d
    \end{cases}.\label{eq:true_soln}
\end{equation}
solves~\eqref{eq:Transmission_PDE} with~$v=v^{\In}$. We may therefore numerically solve \eqref{eq:Transmission_PDE} using the method described in this section and compare the numerical solution to the true solution~\eqref{eq:true_soln}.
In~\figref{fig:single-ana_test}, we show the accuracy when~$\tilde\gamma$ is discretized using 72 16th-order panels and~$L=10$. The maximum error is determined by the adaptive integration tolerance, which is $10^{-12}$. 

\begin{figure}
    \centering
    \includegraphics[width=.6\textwidth]{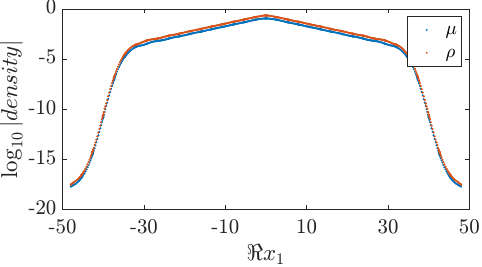}
    \caption{The modulus of the solution densities from the analytic solution test in~\figref{fig:single-ana_test}. The densities from the two halves of~$\tilde\gamma$ overlap because of the symmetry of the problem and solution. The figure clearly shows that contour deformation was enough to ensure that the densities have decayed to machine precision. }
    \label{fig:single_dens}
\end{figure}

To check that the solver does not depend on~$L$ and our truncation distance, we redo the experiment with those distances and the number of points doubled. In this test, the maximum error is still comparable to the adaptive integration tolerance. We also redo the experiments with~$\tilde\gamma$ truncated at~$|\Re x_1| = L+25$. This time the maximum error is~$1.3\times 10^{-7}$ because we did not go far enough out into the complex plane to ensure that the densities had sufficiently decayed. This can easily be seen from the graph of the corresponding densities are shown in~\figref{fig:single_dens}.

\begin{figure}
    \centering
    \includegraphics[width=\textwidth]{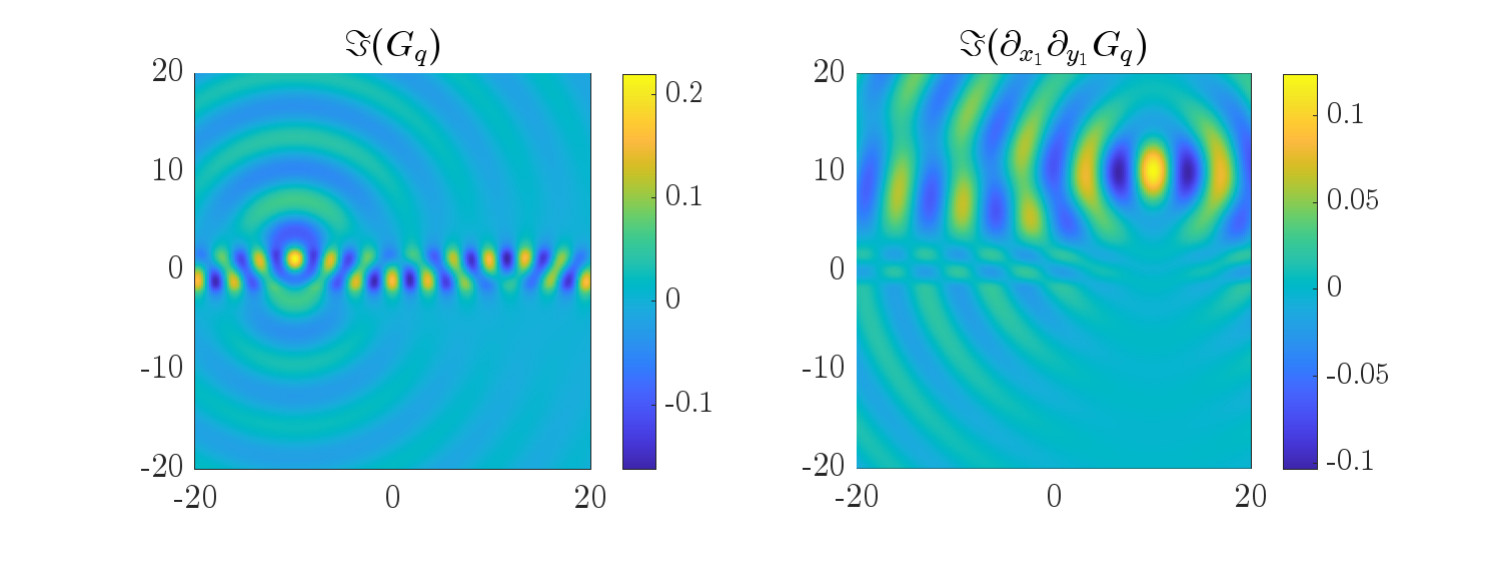}
    \caption{An example of~$G_q$ and one of its derivatives.}
    \label{fig:Gq_examp}
\end{figure}

Now that we have completed the construction of our solver, we are able to evaluate~$G_q$ and it's derivatives. A few examples are shown in~\figref{fig:Gq_examp}.

\section{Solving the matched IE}\label{sec:solve_match}
In this section, we discuss the numerical solution of~\eqref{eqn2.13.001}. In order truncate the problem, we again use coordinate complexification. In particular, we  solve~\eqref{eqn2.13.001} on the deformed contour
\begin{equation}
    \hat\gamma = \{ (0,y_2+i \psi(y_2); 0,x_2+i \psi(x_2))| y_2,x_2\in\bbR\}.
\end{equation}
We defer the proof that~\eqref{eqn2.13.001} is well posed on~$\hat\gamma$ to a future paper. Here we simply observe that the wave-guide Green's function, $G_{q_{l,r}},$ can be analytically continued to~$\hat\gamma$, and so~\eqref{eqn2.13.001} can at least be defined. This follows from the Fourier representation of $G_{q_{l,r}}$ given in~\cite{EpWG2023_1}. 
For the sake of brevity, we defer a rigorous justification to~\cite{CmplxWG2}.

We discretize~\eqref{eqn2.13.001} on~$\hat\gamma$ in the same way we discretized~\eqref{eq:wave_guide_IE} in Section~\ref{sec:wave_discretize}, except that the integral operators are built using the wave-guide Green's functions. We  let~$\bz_i$ be the nodes used to discretize~$\hat\gamma$ and assume that there are~$m$ of them. 

At this point, we have all the tools needed to accurately solve~\eqref{eqn2.13.001}. However, a solver based on this discretization would be slow, as each evaluation of~$G_{q_{l,r}}$ requires us to solve~\eqref{eq:wave_guide_IE}. The two halves of the problem, building the system matrix and recovering the solution, require different approaches to handle this slowness. In the remainder of this section, we explain how those computations can be accelerated,  starting with reconstruction of the solution.

Once we have solved for~$\sigma$ and~$\tau$, the naive evaluation of the field $u_{l}=\cS_{q_l}[\tau]+\cD_{q_l}[\sigma]$ at each target  costs~$O(nm)$. If there~$n_t$ targets, then the cost becomes~$O(n_tnm)$, which is infeasibly large. To accelerate this computation, we  note that~$u_l(\bx)$ is a solution of~\eqref{eq:Transmission_PDE} with
\begin{equation}
    v^{\In}(\bx)=\sum_{i} G_{k_{l,r}(\bx)}(\bx;\bz_i)\tau_i w_i + \partial_{z_1}G_{k_{l,r}(\bx)}(\bx;\bz_i)\sigma_i w_i.
\end{equation}
Evaluating the right hand side of~\eqref{eq:wave_guide_IE} with this incoming field costs~$O(nm)$. Solving~\eqref{eq:wave_guide_IE}  then costs~$O(n)$ and evaluating~\eqref{eq:wave_guide_rep} costs~$O(n_t n)$. If we add the cost of evaluating~$v^{\In}$ at all the targets,~$O(n_tm)$, then we can see that the total cost to evaluate~$u_{l}$ can be brought from~$O(n_tnm)$ down to~$O(n_t(n+m) + nm )$. This lower cost makes it possible to create high resolution plots of the solution. If desired, one could accelerate this to linear complexity using an FMM~\cite{zFMM}. A nice feature of this approach, is that it is completely general. We can therefore use it for any decomposable problem.

We now turn to solving~\eqref{eqn2.13.001}. The simplest approach would be to avoid building the system matrix entirely, and solve~\eqref{eqn2.13.001} using an iterative solver. We could then use the strategy described above to apply most of the system matrix quickly. The singular quadrature for close sources and targets would have to be handled separately, which would greatly slow down the algorithm.

Instead, we  build a compressed representation of the wave-guide Green's function. For~$\bz_i$ far from the wave-guides, this representation is built using low rank factorizations. When those are not available, we use Chebyshev expansions. In the remainder of this section we spell out these steps in more detail. To keep things simple, we build separate compressions for~$G^{\out}_{q_{l}}$ and~$G^{\out}_{q_{r}}$ and suppress the dependence on~$l,r$ in our discussion.

\subsection{Outgoing skeleton}\label{sec:inskel}

When building the system matrix, we must compute
\begin{equation}
    v(\bz_i) = \cS_k[\rho] + \cD_k[\mu]\approx \sum_{j=1}^n G_k(\bz_i;\bx_j)w_j\rho_j +\partial_{\nu_j}G_k(\bz_i;\bx_j) w_j\mu_j
\end{equation}
for every~$i$ and with a different pair of densities for each source. We therefore desire a way to evaluate this sum quickly. 

We let~$E$ be the matrix mapping from~$\rho_j$ and~$\mu_j$ to~$v(\bz_i)$. Let~$Z_F=\{\bz_i\}_{i\in I_F}$ be the set of~$\bz_i's$ that are at least a wavelength from any of the~$\bx_j$ and~$E_F = E(I_F,:)$. The submatrix~$E_F$ is numerically low rank provided the extent of~$\bx_j$'s and~$\bz_i$'s is not too large (see~\cite{zFMM}). We can therefore apply a large slice of~$E$ using a low rank factorization of~$E_F$.

A particularly nice factorization, is the interpolative decomposition
\begin{equation}
    E_F \approx E_F(:,J_{\rm{out-skel}}) T, \label{eq:fac_def}
\end{equation}
where~$J_{\rm{out-skel}}$ represents an ``important" set of sources and~$T$ is an interpolation matrix. We refer to~$J_{\rm{out-skel}}$ as the \textit{outgoing skeleton} of~$\tilde \gamma$. Physically, the factorization \eqref{eq:fac_def} is equivalent to the statement
\begin{equation}
    v(\bz_i) \approx \sum_{j\in J_{\rm{out-skel}}} \cG_j(\bz_i; \bx_j) \left[T\begin{pmatrix}
        \rho\\\mu
    \end{pmatrix}\right]_j,\label{eq:fac_fct_def}
\end{equation}
for~$i\in I_F$, where~$\cG_j$ is~$G_k$ or~$\partial_{\bs n_j}G_k$. 

In practice, this idea has a few drawbacks. Namely that we must build a different factorization for~$\partial_{z_2}v$ and if we change~$\hat\gamma$. It is therefore advantageous to design~$J_{\rm{out-skel}}$ and~$T$ such that~\eqref{eq:fac_fct_def} holds for a more flexible choice of~$\bz$'s. A common strategy when all the points are in~$\bbR^2$ is to ensure that~\eqref{eq:fac_fct_def} holds for all~$\bz$ on a curve~$\gamma_P$, called a proxy curve, that encompasses the points in~$Z_F$ (see~\cite{martinsson2019fast}). A fairly simple argument using Green's identity shows that if~\eqref{eq:fac_fct_def} holds on~$\gamma_P$, then it holds for every point inside~$\gamma_P$ and also for the derivatives of~$v$. We shall use use this idea to build our factorization, with~$\gamma_P$ being a curve that encompasses the real parts of the points in~$Z_F$ and is at least a wavelength from any~$\bx_j$. 

The proof that a skeleton built using this~$\gamma_P$ gives correct answers for complex~$\bx$ and~$\bz$ is the subject of ongoing work. We therefore present numerical evidence that the proxy curve is applicable to this case. We choose~$\gamma_P$ to be the union of boundaries of the boxes
$$B_+=[-1,1]\times [d+2\pi/k,100]\text{ and }B_-=[-1,1]\times[-100,-d-2\pi/k].$$
We build the outgoing skeleton such that~\eqref{eq:fac_fct_def} holds to a tolerance~$\epsilon$ for 300 points on each segment of~$\gamma_P$.

\begin{figure}
    \centering
    \includegraphics[width=0.6\textwidth]{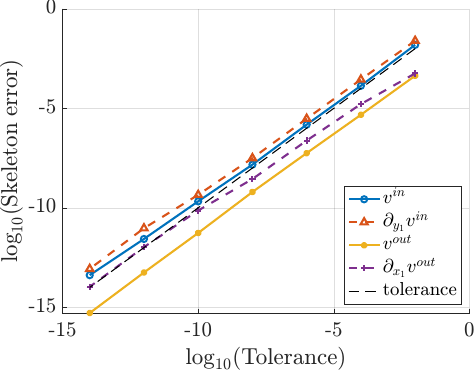}
    \caption{The error in the formulas~\eqref{eq:fac_fct_def} and~\eqref{eq:in-skel} as a function of the tolerance used to build the skeletons. The error presented is a maximum over points spread out in the boxes~$B_\pm$.}
    \label{fig:skel_err}
\end{figure}

We then test the skeletonization error, i.e. the error in~\eqref{eq:fac_fct_def}, on a grid of 250 targets in each~$B_{\pm}$. \figref{fig:skel_err} shows how the maximum interpolation error depends on the tolerance used to build the skeleton. As we can see, the interpolation error matches the tolerance, and so the proxy curve can be used with complex~$\bx_j$. Similar results are found when we take the 
$x_1$-derivative of~\eqref{eq:fac_fct_def}, which will be required to evaluate~$\partial_{x_1}v$. We also added imaginary parts~$\Im z_2 = \psi(\Re z_2)$ to the targets and found identical results.

By refining~$\tilde\gamma$, we also observed the the number of skeleton points is independent (within~$\pm1$) of the number of panels in~$\tilde\gamma$.

\subsection{Incoming skeletons}\label{sec:outskel}

Building the system matrix also requires us to evaluate many different right hand sides for~\eqref{eq:wave_guide_IE}. Namely, we must evaluate
\begin{equation}
    f_{ij} :=  G_k(\bx_i;\bz_j) -G_{k_1}(\bx_i;\bz_j)
\end{equation}
and
\begin{equation}
    g_{ij} := \partial_{\bs n_i}G_k(\bx_i;\bz_j) -\partial_{\bs n_i}G_{k_1}(\bx_i;\bz_j)
\end{equation}
for every~$i$ and~$j$. We denote the combined matrix of~$f_{ij}$ and~$g_{ij}$'s by~$R$. The slice of~$R$ corresponding to~$\bz_j\in Z_F$ is~$R_F$. This slice  has a low rank factorization for the same reason as~$E_F$ does. We shall therefore work to build its interpolative decompositions
\begin{equation}
R_F \approx X R(I_{\inskel},:).
\end{equation}
We refer to~$I_{\inskel}$ as the \textit{incoming skeleton} of~$\tilde \gamma$. Since~$G_k-G_{k_1}$ does not satisfy a Green's theorem, we can not directly use a proxy curve to build this factorization. Instead, we therefore build a skeleton that works for each kernel separately, i.e. we choose~$I_{\inskel}$ and~$X$ such that
\begin{equation}
    G_{\tilde k}(\bx_i;\bz_j) \approx \sum_{l\in I_{\inskel}} X_{il} \tilde\cG_l(\bx_l;\bz_j)  \label{eq:in-skel}
\end{equation}
and
\begin{equation}
    \partial_{\bs n_i}G_{\tilde k}(\bx_i;\bz_j) \approx \sum_{l\in I_{\inskel}} X_{(i+n)l} \tilde\cG_l(\bx_l;\bz_j),
\end{equation}
where~$\tilde\cG_l=G_{\tilde k}$ or~$\partial_{\bs n_l}G_{\tilde k}$ and~$\tilde k=k,k_1$. Since~$G_k$ and~$G_{k_1}$ independently satisfy Green's identities, we can build this factorization using the same proxy curve~$\gamma_P$ as we used to build the outgoing skeleton.

\begin{remark}
    The incoming skeleton described here could be used an outgoing skeleton. In practice, we chose not to reuse the skeletons to keep the outgoing skeleton as small as possible.
\end{remark}

We test our incoming skeleton, in the same way as the outgoing section in the previous section. \figref{fig:skel_err} shows that the skeletonization error matches the tolerance used to compare the skeleton, as we would hope. Finally, we observe that the same skeleton can be used to compute~$\partial_{y_1} G_q^{\out}$.

\subsection{Far sources and far targets}
In this section, we show how the incoming and outgoing skeletons can be used to build a fast evaluator for~$G_q^{\out}(\bx;\by)$ when~$\bx,\by\in Z_f$. As described in the previous section, for any source~$\by$ in~$B$, we can efficiently evaluate the right hand side of~\eqref{eqn2.13.001} by applying~$X$ to the vector
\begin{equation}
    \rhs_j(\by) = [[\tilde\cG_j(\bx_j;\by)]]_{\tilde\gamma}.   
\end{equation}
We can therefore compute the densities as
\begin{equation}
     (\mu(\bz);\rho(\by))= K^{-1}X \rhs(\by),\label{eq:fast_density}
\end{equation}
where~$K$ is the system matrix for~\eqref{eq:wave_guide_IE}.
This formula is helpful because~$K^{-1}X$ can be precomputed allowing us to efficiently solve~\eqref{eq:wave_guide_IE} for any source~$\by\in Z_f$. We can efficiently evaluate the potentials generated by these densities using~\eqref{eq:fac_fct_def}. Putting it all together we have the formula
\begin{equation}
G_q^{\out}(\bx;\by) \approx \sum_{j\in J_{\rm{out-skel}}} \cG_j(\bx; \bx_j) \left[(TK^{-1}X) \rhs(\by)\right]_{j},\label{eq:far2far}
\end{equation}
which is accurate for any~$\bx,\by\in Z_F$. If we precompute~$TK^{-1}X$, then this formula can be used to evaluate~$G_q^{\out}$ in~$O(n_{\inskel}n_{\rm{out-skel}})$ time, independent of the number of points in~$\tilde \gamma$.

\begin{figure}
    \centering
    \includegraphics[width=0.6\textwidth]{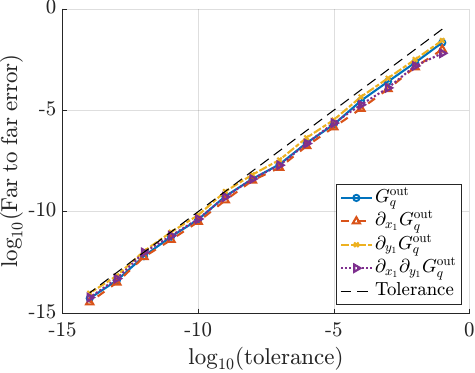}
    \caption{This figure shows how the error when~\eqref{eq:far2far} is used to approximate~$G_q^{\out}(\bx;\by)$ and its derivatives depends on the skeletonization tolerance. The error reported is the maximum pointwise error over points spread out in the boxes~$B_\pm$.}
    \label{fig:far2far_err}
\end{figure}

We test this approach by finding the error in computing~$G_q^{\out}(\bx;\by)$ for every pair of points used in Section~\ref{sec:inskel}. \figref{fig:far2far_err} shows how the maximum of these errors scales with the tolerance used to compute the skeletons. To further test this method, we check the error in the equivalent formulas for~$\partial_{x_1}G_q^{\out}, \partial_{y_1}G_q^{\out}$, and~$\partial_{x_1}\partial_{y_1}G_q^{\out}$. We can see that all of the errors scale with the skeletonization tolerance correctly.

\subsection{Far sources and near targets}
When the target~$\bx$ is not in~$Z_F$, we are not able to use the outgoing skeleton to compute~$G_q^{\out}(\bx;\by)$. We therefore need a different way to quickly evaluate~$G_q^{\out}$. When~$\by\in Z_F$, \eqref{eq:fast_density} tells us that~$G_q^{\out}(\bx;\by)$ can be written as a linear combination of the potentials generated by the densities that are the columns of~$K^{-1}X$. It therefore remains to find an efficient way to generate those potentials.

Let~$v_i$ be the field generated by the~$i$th column of~$K^{-1}X$. In order to rapidly evaluate those fields, we shall generate piecewise Chebyshev interpolants of~$u_i(0,\cdot)$. Let~$I_1,\ldots,~I_N$ be a collection of intervals that partition the interval~$[-d-2\pi/k,d+2\pi/k]$ such that~$q$ is constant on each interval. Typically, we shall use one interval per region that~$q$ is constant, but wide or high frequency wave-guides may need more intervals.

On each of the intervals~$I_1,\ldots,I_N$ we build a Chebyshev expansion. If~$l_j$ and~$c_j$ are length and center of~$I_j$, then this takes the form.
\begin{equation}
     v_i(0,x_2)\approx \sum_{l=1}^p a_{i,l,j} T_l\lp2\frac{x_2-c_j}{l_j}\rp \quad\forall x_2\in I_j.\label{eq:cheb_far2near}
\end{equation}

In general, we have no guarantee that the~$v_i$'s are nice functions in the neighbourhood of~$\tilde\gamma$. We therefore do not expect $v_i(0,\cdot)$ to be smooth so  these Chebyshev expansions may not converge. Fortunately, we do not need to evaluate individual~$v_i$'s. Instead we only use linear combinations of the~$v_i$'s. Indeed~\eqref{eq:fast_density} is equivalent to the statement that
\begin{equation}
    G_q^{\out}(0,x_2;\by)\approx \sum_{i\in I_{\inskel}} v_i(0,x_2)\cG_i(\bx_i;\by)
\end{equation}
for all~$x_2\in \hat\gamma$. If we plug our expansions of~$v_i$ into this expression, then we find\begin{equation}
    G_q^{\out}(0,x_2;\by)\approx \sum_{l=1}^p\lp\sum_{i\in I_{\inskel}} a_{i,l,j}\cG_i(\bx_i;\by) \rp T_l\lp2\frac{x_2-c_j}{l_j}\rp\label{eq:cheb_far2near_G}
\end{equation}
for all~$x_2\in I_j$. Equation~\eqref{eq:cheb_far2near_G} is a Chebyshev interpolant of~$G_q^{\out}(0,x_2;\by)$, up to the tolerance used in the skeletonization and solving~\eqref{eq:wave_guide_IE}. Since~$\by\in Z_F$, the function~$G_q^{\out}(0,x_2;\by)$ is a solution of~\eqref{eq:Transmission_PDE} with smooth boundary data. The Chebyshev interpolant~\eqref{eq:cheb_far2near_G} therefore converges spectrally, making~\eqref{eq:cheb_far2near_G} an efficient and accurate method for computing~$G_q^{\out}(0,x_2;\by)$.

\begin{figure}
    \centering
\includegraphics[width=0.5\textwidth]{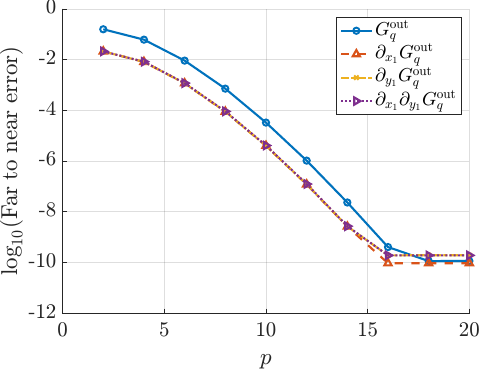}
    \caption{This figure shows how the error in~\eqref{eq:cheb_far2near} depends on the order of the Chebyshev expansions.}
    \label{fig:far_2_near}
\end{figure}

To test this algorithm, we build the incoming skeletonization with a tolerance of~$10^{-10}$ and build the Chebyshev expansions~\eqref{eq:cheb_far2near} with various lengths~$p$. We test the interpolation accuracy using the same grid of sources used to test the incoming skeleton and~$100$ targets evenly spaced on the line $\{0\}\times (\cup_j I_j)$. The maximum errors are reported in~\figref{fig:far_2_near}. We can see that the errors do indeed converge spectrally with~$p$ up to the skeletonization tolerance.

Since some of our matrix entries need~$x_1$ derivatives of~$G_q$, we also build Chebyshev expansions of~$\partial_{x_1}v_i(0,x_2)$. We then test the equivalent formulas to~\eqref{eq:cheb_far2near_G} for all of the required derivatives of~$G_q^{\out}$. \figref{fig:far_2_near} shows that these interpolants are also spectrally accurate.

In theory, the interpolation error can also be determined by the discretization error in~\eqref{eq:wave_guide_IE}. However, since~$\by$ is far from~$\tilde \gamma$, the densities are well resolved and adaptive integration is used to ensure that the~$v_i$'s are computed accurately.

\subsection{Near sources and far targets}
We can use a similar idea to the previous section when~$\bx\in Z_F$ but~$\by$ is not. By~\eqref{eq:fac_fct_def}, for any~$\bx\in Z_F$ we have that
\begin{equation}
    G_q^{\out}(\bx;0,y_2) \approx \sum_{i\in J_{\rm{out-skel}}} \cG_i(\bx; \bx_i) \left[T\begin{pmatrix}
        \rho (0,y_2)\\\mu (0,y_2)
\end{pmatrix}\right]_i\label{eq:near_2_far_setup}
\end{equation}
We build piecewise Chebyshev expansions of the density skeletons on each subinterval:
\begin{equation}
    \left[T\begin{pmatrix}
        \rho (0,y_2)\\\mu (0,y_2)
    \end{pmatrix}\right]_i = \sum_{l=1}^p b_{i,l,j}T_l\lp2\frac{y_2-c_j}{l_j}\rp \quad\forall y_2\in I_j.\label{eq:cheb_near2far}
\end{equation}
Plugging these expansions into~\eqref{eq:near_2_far_setup} gives
\begin{equation}
    G_q^{\out}(\bx;0,y_2) \approx \sum_{l=1}^p \lp \sum_{i\in J_{\rm{out-skel}}} b_{i,l,j}\cG_i(\bx; \bx_i) \rp T_l\lp2\frac{y_2-c_j}{l_j}\rp, \label{eq:near_2_far}
\end{equation}
which is a Chebyshev expansion on the interval~$I_j$. Since~$G_q^{\out}(\bx;0,y_2)$ is a piecewise smooth function of~$y_2$, this expansion  converges spectrally up to the error in the skeletonization and the discretization of~\eqref{eq:wave_guide_IE}.

\begin{figure}
    \centering
    \begin{subfigure}{0.45\linewidth}
    \includegraphics[width=0.9\textwidth]{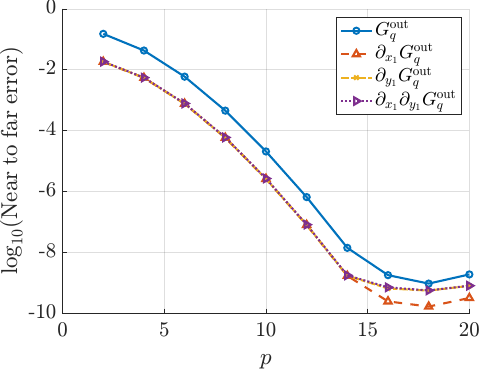}
    \caption{Maximum near to far interpolation error as a function of the order of the Chebyshev expansions for the various derivatives of~$G_q^{\out}$}
    \label{fig:near_2_far}
    \end{subfigure}\hspace{.1cm}
    \begin{subfigure}{0.45\linewidth}
    \includegraphics[width=0.9\textwidth]{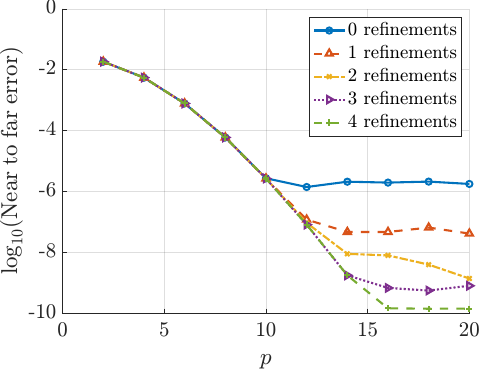}
    \caption{Maximum near to far interpolation error as a function of the order of the Chebyshev expansions for~$\partial_{x_1}\partial_{y_1} G_q^{\out}$ as a function of the number of refinements of the boundary~$\tilde\gamma$.}
    \label{fig:near_2_far_conv}
    \end{subfigure}
    \caption{These figures shows how the error in~\eqref{eq:near_2_far} various with the expansion order and the discretization of~$\tilde\gamma$.}
    \label{fig:near2far_tot}
\end{figure}

To test this algorithm, we build the outgoing skeletonization with a tolerance of~$\epsilon=10^{-10}$ and build the Chebyshev expansions~\eqref{eq:cheb_near2far} with various lengths~$p$. We compute the interpolation accuracy using the same grid of targets used to test the outgoing skeleton and~$100$ sources evenly spaced on the line $\{0\}\times (\cup_j I_j)$. The maximum errors are reported in~\figref{fig:near_2_far}. We can see that the errors do indeed converge super-algebraically with~$p$.

We can use similar formulas to compute the derivatives~$G_q^{\out}$, except that new expansions must be build for the~$y_1$ derivative of the densities. We test the accuracy of these for all of the required derivatives of~$G_q^{\out}$. \figref{fig:far_2_near} shows these interpolants also work as expected.

In all of these experiments, the error stagnates at a value higher than would caused by the skeletonization. This is due to the densities not being sufficiently resolved when the sources and Chebyshev nodes are quite close to~$\tilde\gamma$. To verify this observation, we repeat the same experiment with several different discretizations of~$\tilde\gamma$ and the highest order required derivative of~$G_q^{\out}$. \figref{fig:near_2_far_conv} shows that the stagnation point does indeed depend on the discretization and can be pushed down to the skeletonization tolerance. When we solve~\eqref{eqn2.13.001} in Section~\ref{sec:scatter_mat}, this experiment  helps us to decide how much to refine~$\tilde\gamma$.

\begin{remark}
    It is sometimes the case that we would like to know the solution a great distance from the waveguide. We can achieve this by iteratively extending the computational domain in a manner similar to that described in~\cite{hoskins2020discretization}. This is cheaper than solving directly, as it only requires us to evaluate the waveguide Green's functions at targets far from the waveguides. Due to the way we have constructed the skeletons, they can be reused in each iteration.
\end{remark}

\subsection{Near sources and near targets}\label{sec:near2near}
We now turn to the final case. When neither~$\bx$ nor~$\by$ are in~$Z_F$, we cannot use skeletonization to build efficient evaluators. We therefore simply use tensor product Chebyshev expansions for each pair of subintervals~$I_i$ and~$I_j$. We show below that there are more efficient compression schemes, but this is sufficient for our purposes. Since we only build these expansions for the region close to~$\tilde\gamma$, we do not need too many subintervals and the precomputation is not  too costly. It remains, however, to analyze their accuracy.

The accuracy is determined by the smoothness of $G_q^{\out}(0,x_2;0,y_2).$ The following lemma summarizes some results in~\cite{EpWG2023_1}.
\begin{lemma}[\cite{EpWG2023_1}]\label{lem:Green_smooth}
The function~$G_q^{\out}(0,x_2;0,y_2)$ and it's derivatives are smooth in each region where~$q$ is constant apart from a singularity as~$x_2,y_2\to \pm d$.

In that limit the derivatives of~$G_q^{\out}$ satisfy
\begin{equation}
    \partial_{\bx,\by}^\alpha G_q^{\out}(0,x_2;0,y_2)\sim (x_d+y_d)^{2-|\alpha|}\log \lp x_d+y_d \rp ,
\end{equation}
where~$\alpha$ is a multi index of derivatives and 
\begin{equation}
    x_d = ||x_2|-d| \quad \text{and} \quad y_d = ||y_2|-d|.
\end{equation}
\end{lemma}
The statement in Lemma~\ref{lem:Green_smooth} is slightly different than is given in~\cite{EpWG2023_1} because we have chosen $v^{\inc}=G_{k(\bx)}$ rather than $G_k$.

\begin{figure}
    \centering
    \includegraphics[width=0.6\textwidth]{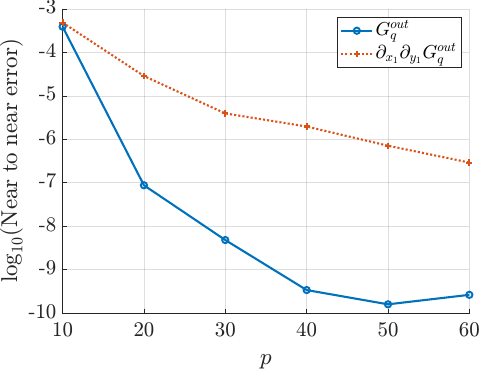}
    \caption{Our near to near test}
    \label{fig:near_to_near}
\end{figure}

To test our compressed representation, we place 100 equally spaced points on the line~$\{0\}\times [-d-2\pi/k,d+2\pi/k]$ and compute~$G_q^{\out}(0,x_2;0,y_2)$ for every pair of points. Due to the singularities indicated by Lemma~\ref{lem:Green_smooth}, the Chebyshev expansions do not converge in an~$L^\infty$ sense. We therefore report an approximation of the $L^2$-error in \figref{fig:near_to_near}. We see that the error converges to expected value as we increase~$p$, albeit slowly. We could obtain better convergence by building the singularities in Lemma~\ref{lem:Green_smooth} into our expansion. However, the errors shown in \figref{fig:near_to_near} are small enough for our purposes.

\section{Numerical experiments}\label{sec:test_match}
We now test the numerical solver developed in the previous section. We choose~$k= 1$ and set $k_{1;l}=2, d_l=2; k_{1;r}=3, d_r=4$ so that
\begin{equation}
    q_l(x_2) = (2^2-1) \chi_{[-2,2]}(x_2)\qand q_r(x_2) =  (3^2-1) \chi_{[-4,4]}(x_2).
\end{equation}
We discretize~$\hat\gamma$ using 96 16th order panels and~$L=10$. We build the skeletonization with a tolerance of~$10^{-10}$ and build the expansions~\eqref{eq:cheb_far2near} and~\eqref{eq:cheb_near2far} to order 20. The other expansions are built to  order 60.

\begin{figure}
    \centering
    \begin{subfigure}{0.45\linewidth}
        \centering
    \includegraphics[width=\textwidth]{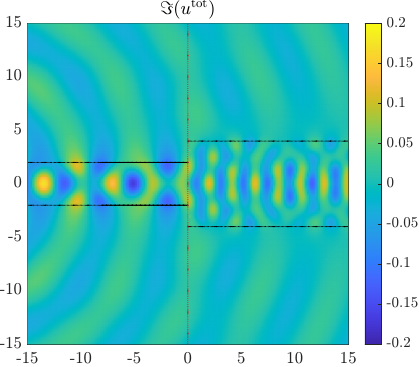}
    \caption{The imaginary part of the solution~$u$.}
    \label{fig:match_soln}
\end{subfigure}\hspace{.1cm}
    \begin{subfigure}{0.45\linewidth}
    \centering
    \includegraphics[width=\textwidth]{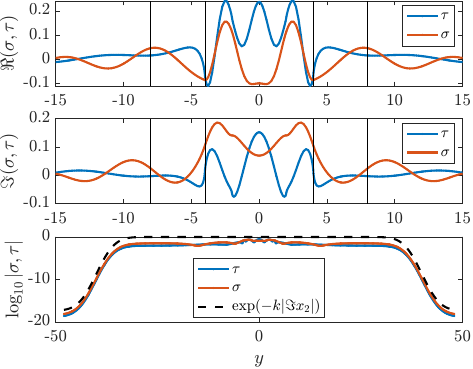}
    \caption{The corresponding densities~$\sigma$ and~$\tau$.}
    \label{fig:match_density}
    \end{subfigure}
    \caption{These figures show a solution of the matched wave-guide problem in the case where two wave-guides meet along the~$x_2$-axis.}
\end{figure}

\begin{figure}
    \centering
    \includegraphics[width= 8cm]{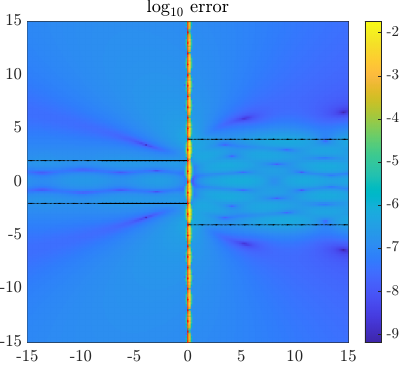}
    \caption{The estimated error in our self convergence study. We can see that the solution is accurate to 6 digits everywhere away from~$\hat\gamma$.}
    \label{fig:match_self_conv}
\end{figure}

 We set
 \begin{equation}
    v^{\In}_l(\bx) = -G_{q_l}(\bx;\bx_0),
\end{equation}
where~$\bx_0=(-10,0)$, and~$v^{\In}_r=0$. Figure~\ref{fig:match_soln} shows the resulting total field. To test the accuracy of our solve, we re-solve this problem with twice as many points on the contour~$\hat\gamma$. The error in this experiment is shown in~\figref{fig:match_self_conv}. As we can see, the error is quite small away from~$\hat\gamma$. The solution is less accurate around~$\hat\gamma$ because the integrals in~\eqref{eq:u_rep} become nearly singular, something that could be resolved using adaptive integration. We thus conclude that we have solved our problem to the accuracy allowed by our compression algorithms.

It is interesting to observe the behaviour of the densities~$\sigma$ and~$\mu$. From~\figref{fig:match_density}, we can make two observations. The first is that the contour deformation has indeed caused the densities to exponentially decay and so we have truncated~$\hat\gamma$ in a reasonable location. The second observation is that the densities are smooth on each side of boundaries of the wave-guides. This is a nice feature of our integral equation. Even though the underlying PDE has a corner at those points, the integral equation is not singular there.

To further test our discretization, we use an analytic solution test similar to the one in Section~\ref{sec:test_single}. We redo the experiment with~$\bx_0=(10,0)$. In this case the solution of~\eqref{eqn1.2} is
\begin{equation}
    u(\bx) = \begin{cases}
        G_{q_l}(\bx;\bx_0) & \text{ if }x_1<0\\
        0 & \text{ if }x_1>0.
    \end{cases}
\end{equation}
The errors in the numerical solution are comparable to the self-convergence error estimate in~\figref{fig:match_self_conv}, except they are about digit better. We can thus conclude that we are indeed able to accurately solve~\eqref{eqn1.2}. The increased accuracy is due the fact that the solution, in this, case does not have a singularity at the intersection of~$\tilde\gamma$ and~$\hat\gamma$, whereas the solution in the previous test has weak singularities there. 

To verify that the choice of~$L$ in~$\psi$ and our choice of
truncation do not affect the accuracy, we redo the analytic solution
test with~$L$, the truncation distance, and the number of points
doubled. There error is still bounded by $10^{-6}$ and so we have
chosen these parameters sufficiently large.

\section{The scattering matrix and radiated power} \label{sec:scatter_mat}

In the design process, it is often important to know how much power is transferred through a junction,  how much is reflected, and how much is lost to radiation. This information is encoded in the scattering matrix. With the ability to compute scattering matrices, we can build a very efficient simulation of an optical circuit. In this section, we describe how to build the scattering matrix and how to use it to measure the radiated power.

\subsection{Finding the wave-guide modes}
As described above, the wave-guide modes are of the form
\begin{equation}
    u_{j;l,r,\pm} (x_1,x_2) = e^{\pm i\xi_{j;l,r}} v_{j;l,r}(x_2).
\end{equation}
where the~$v_{j;l,r}(x_2)$'s are the solution of~\eqref{eqn2.2.002}. The~$\xi_{j;l,r}$'s are the roots of~\eqref{eq:mode_eqn} that lie in the interior of~$(k,k_{1;l,r})$.   To find these roots, we use the smoother equivalent form of equation~\eqref{eq:mode_eqn}
\begin{multline}\label{eq:smth_mode}
    f_{l,r}(\xi) = \lp k_{1;l,r}^2+k^2-2\xi^2\rp  \sin\lp 2d\sqrt{k_{1;l,r}^2 - \xi^2} \rp \\
    - 2\sqrt{\xi^2-k^2}\sqrt{k_{1;l,r}^2 - \xi^2}\cos\lp2d\sqrt{k_{1;l,r}^2 - \xi^2} \rp=0.
\end{multline}
Not that while $k_{1;l,r}$ is always a root of this equation it is never a root of the Wronskian.

To find the roots of~$f_{l,r}$, we build a 500 term Chebyshev expansion of each function. To avoid a loss of accuracy caused by the square root singularities, we build this expansion on the interval~$(k+0.01,k_{1;l,r}-0.01)$, which allows us to find all of the modes for most wave-guides. Wave-guides with modes outside this interval could be handled with a longer expansion over a bigger interval. We find the roots of those expansions using their colleague matrices (see~\cite{good1961colleague}). Any roots of those expansions that lie in~$(k+0.01,k_{1;l,r}-0.01)$ are the frequencies of the wave-guide modes. 

Now that we have found these modes, we can compute the projection of~$u$ onto them. The projection onto~$u_{j;r,+}$ is defined to be
\begin{equation}
    c^r_j = \int_{-\infty}^\infty \overline{u_{j;r,+}}(x_1,x_2) u(x_1,x_2) dx_2 \label{eq:proj_def}
\end{equation}
for any~$x_1>0$. The projections onto~$u_{j;l,-}$ are defined similarly. 
By equation (206) in~\cite{EpWG2023_1}, these coefficients can also be computed by the formula
\begin{equation}
    c^r_j = \int_{-\infty}^\infty v_{r,j}(x_2) \lp -i\xi^r_j \sigma +  \tau \rp dx_2.\label{eq:proj_form}
\end{equation}
As further verification of our solver, we compute the projections of~$u$ in~\figref{fig:match_soln} using both formulas. We discretize~\eqref{eq:proj_def} using~$x_1=5$ and~$x_1=25$ and eight 16th order Gauss-Legendre panels on the numerical support of the wave-guide modes. The the difference in each of the computed coefficients was less than~$10^{-10}$.


\subsection{Energy calculations}\label{sec:NRG_calc}
\begin{figure}
    \centering
    \includegraphics[width= 12cm]{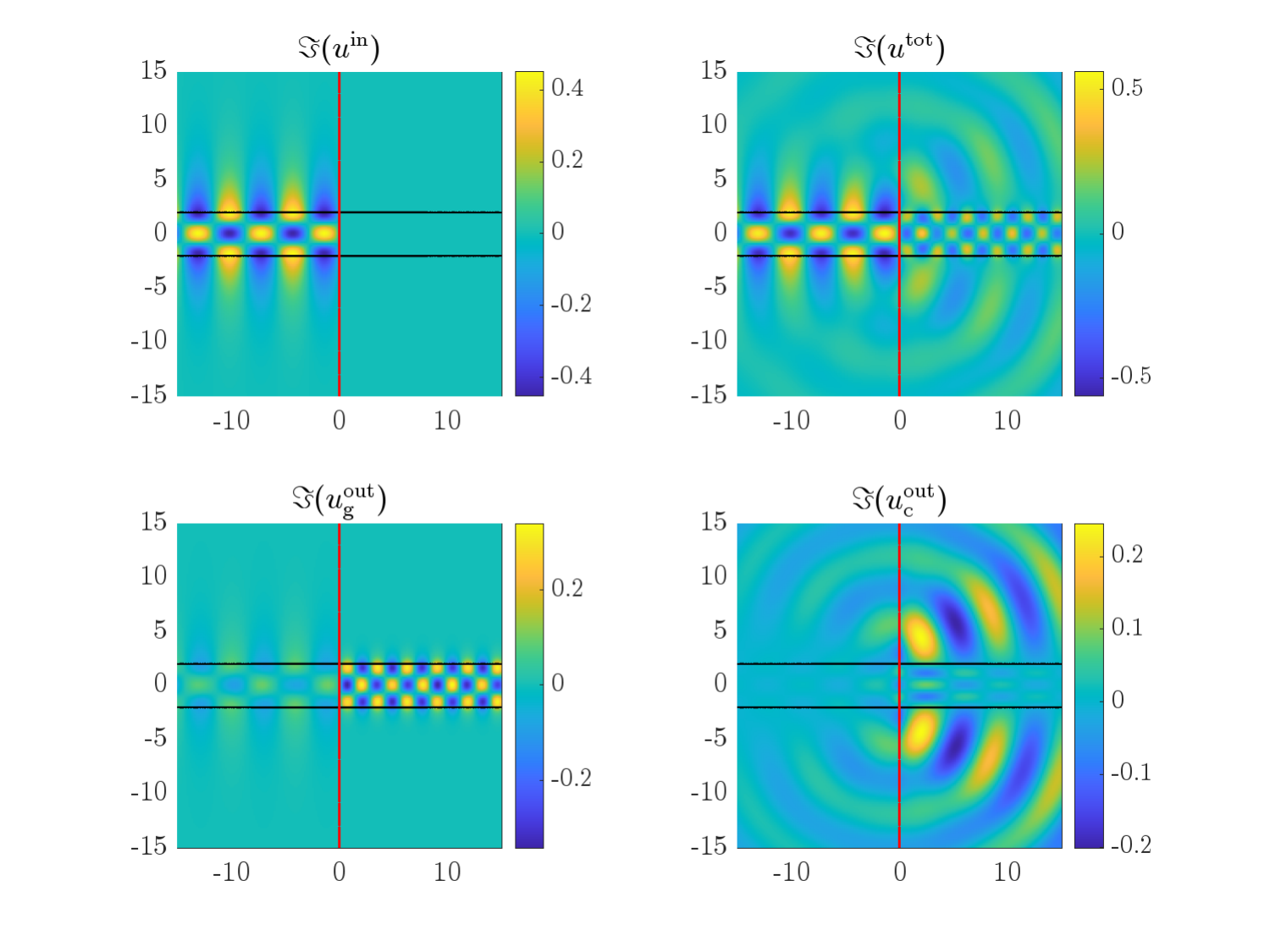}
    \caption{In this figure, we demonstrate how the solution can be decomposed into guided-mode and radiated components. The top two figures show the incoming and total fields. The bottom left figure shows the guided-mode component~$u_g^{\out}$. The bottom right shows the radiated part~$u_c^{\out}$ computed as~$u^{\out}-u_g^{\out}$.}
    \label{fig:radiating_split}
\end{figure}

Now that we have  computed the wave-guide modes, we can solve~\eqref{eqn2.13.001} with~$u^{\In}$ equal to~$u_{j;l,+}$. \figref{fig:radiating_split} shows the resulting~$u$ and the  total field. Using~\eqref{eq:proj_form}, we can compute the guided-mode components of~$u^{\out}$ by
\begin{equation}
    u^{\out}_g = \lp\sum_j c^l_{j,-}u_{j;l,-}\rp \chi_{x_1<0}+ \lp\sum_j c^r_{j,+}u_{j;r,+}\rp\chi_{x_1>0} .
\end{equation}
The difference~$u^{\out}_c:=u^{\out}-u^{\out}_g$ represents the portion of the field that is radiated away from the junction. The components~$u^{\out}_g$ and~$u^{\out}_c$ are shown in~\figref{fig:radiating_split}. This calculation shows that a much larger portion of the wave is transmitted than reflected. We can also see that a large portion is radiated away. This calculation represents the information in one column of the scattering matrix introduced in Section~\ref{ss.scat_matrix.003}. 

With these projections, we are able to determine how much power is sent into each wave-guide and how much is lost. Using the formulas from Section~\ref{ss.scat_matrix.003} we see that the transmitted power is~$P_r:=\sum_j \xi_j^l |c^r_{j,+}|^2$. The total reflected power is~$P_l:=\sum_j \xi_j^l |c^l_{j,-}|^2$. The total outgoing power is
\begin{equation}
   P^{\out}= P_r+P_l+P_{\rad}:= \lim_{R\to\infty}\Im\int_{|x|=R} \overline{u^{\out}} \partial_{r} u^{\out} \,Rd\theta\label{eq:rad_p_def}
\end{equation}
$\bs \nu$ is the outward normal to~$\Gamma$. Hence the power lost to radiation is
\begin{equation}
 P_{\rad}=P^{\out}-P_r-P_l.
\end{equation}

This definition is slow to evaluate. 
If the incoming data is a single wave-guide mode, $e^{i\xi_j^l x_1}v_j(x_2),$ then Theorem~\ref{scat_unit_thm} gives that
\begin{equation}
    \xi_j^l = P_r + P_l + P_{\rad}. \label{eq:power_bal}
\end{equation}
This formula gives us a much easier way to compute the radiated power, since~$P_r$ and~$P_l$ are easy to compute using~\eqref{eq:proj_form}. To verify that our computed solution satisfies this conservation law, we compute the integral \eqref{eq:rad_p_def} on the boundary of the box centered at the origin with side length 40. We discretize the integral using 500 Gauss-Legendre nodes on each side and found that the computed results satisfy~\eqref{eq:power_bal} to 5 digits. To test that our solver satisfies~\eqref{eq:rad_p_def}, we compute~$P_{\rad}$ using an equivalent integral and find
\begin{equation}
    P_{l}/\xi_j^l = 0.03772, \quad P_r/\xi_j^l = 0.56841,\qand  P_{\rad}/\xi_j^l =0.39388.
\end{equation}
Since $ P_r + P_l + P_{\rad}=1.00001 \xi_j^l,$ we see that our solver satisfies~\eqref{eq:power_bal} to 5 digits.

\section{Other examples}\label{sec:other_ex}
So far in this work, we have  explored the solution of a single decomposable problem. In this section we explore a few other examples. The methods are presented without proof, but we do test each method with an analytic solution test analogous to the one performed in Section~\ref{sec:test_match}.
In these examples the exterior wave number $k=1.$
\subsection{Terminating wave-guides}\label{sec:leaky_end}
\begin{figure}
    \centering
    \begin{subfigure}{0.45\linewidth}
    \centering
    \includegraphics[width= .9\textwidth]{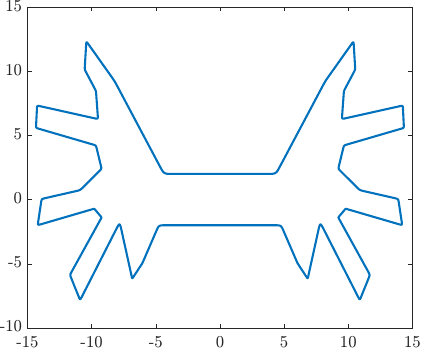}
    \caption{The region~$D$ that we use to terminate a wave-guide. We ensure that~$D$ is symmetric to reduce the number of terms in~\eqref{eqn2.13.001}.}
    \label{fig:cap_shape}
\end{subfigure}\hspace{.15cm}
    \begin{subfigure}{0.45\linewidth}
    \centering
    \includegraphics[width= .9\textwidth]{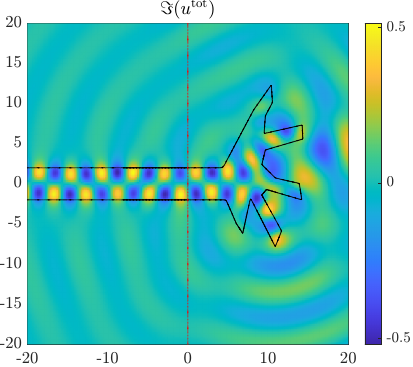}
    \caption{The resulting solution when~$u^{\In}$ is chosen to be a wave-guide mode of the left subdomain.}
    \label{fig:Leaky_tree}
    \end{subfigure}
    \caption{This figure shows our simulation of a terminating wave-guide.}
    \label{fig:tree_examp}
\end{figure}

It is interesting to study how a wave scatters off of the end of a semi-infinite wave-guide. We can model this in the same framework as the matched wave-guide problem if we choose
\begin{equation}
    q_r(\bx) = (2^2-1) \chi_D(\bx),
\end{equation}
where~$D$ is a bounded region that determines the shape of the end of the wave-guide. The Green's function with this choice of~$q_r$ can be computed using the method described in Section~\ref{sec:single_guide} with~$\gamma=\partial D$. As an example, we choose~$D$ to be the region in \figref{fig:cap_shape}. By choosing~$D$ to be symmetric, we ensure that~$\cS_{q_r}'$ and~$\cD_{q_r}$ are still zero, which reduces the computational cost. We solve the resulting modified version of~\eqref{eqn2.13.001} to find the field generated by  a wave-guide mode coming in from the left. The resulting total field is shown in~\figref{fig:Leaky_tree}. An analytic solution test shows that the solution is accurate to 6 digits everywhere away from~$\hat\gamma$ and the method described in Section~\ref{sec:NRG_calc} shows that the solution satisfies~\eqref{eq:power_bal} to 6 digits.

\subsection{Bent Wave-guides}
The algorithm described in previous sections can also be used to study the size of the field radiated when a corner is introduced in a wave-guide. To do this, we choose
\begin{equation}
    q_r(\bx)= (2^2-1) \chi_D(\bx),
\end{equation}
where~$D$ is a wave-guide tilted away from the positive $x_1$-axis by~$\pi/6$ radians, with width chosen to line up with the left wave-guide. The rotated wave-guide Green's function~$G_{q_r}$ can be computed in the same way as the other wave-guides. It should be noted, however, that as the wave-guide is rotated closer to the~$x_2$-axis, the far-field region~$Z_F$ is pushed farther along~$\hat\gamma$ and so the computation becomes more costly.

\begin{figure}
    \centering
    \includegraphics[width= 8cm]{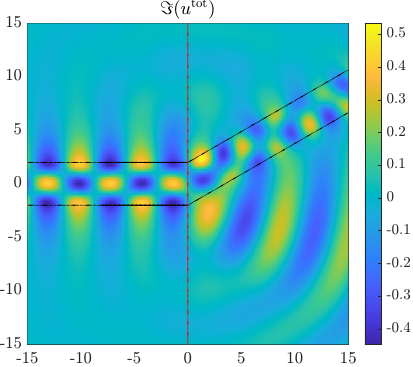}
    \caption{Our simulation of a bent wave-guide with~$u^{\In}$ chosen to be a wave-guide mode.}
    \label{fig:bent_guide}
\end{figure}

With this code, we computed the field generated by an incoming wave-guide mode. The resulting total field is shown in \figref{fig:bent_guide}. An analytic solution tests gives that the solution is accurate to 4 digits everywhere away from the~$\hat\gamma$ and the method described in Section~\ref{sec:NRG_calc} shows that the solution satisfies~\eqref{eq:power_bal} to 5 digits.

\subsection{A Wave-guide junction}

\begin{figure}
    \centering
    \includegraphics[width= 8cm]{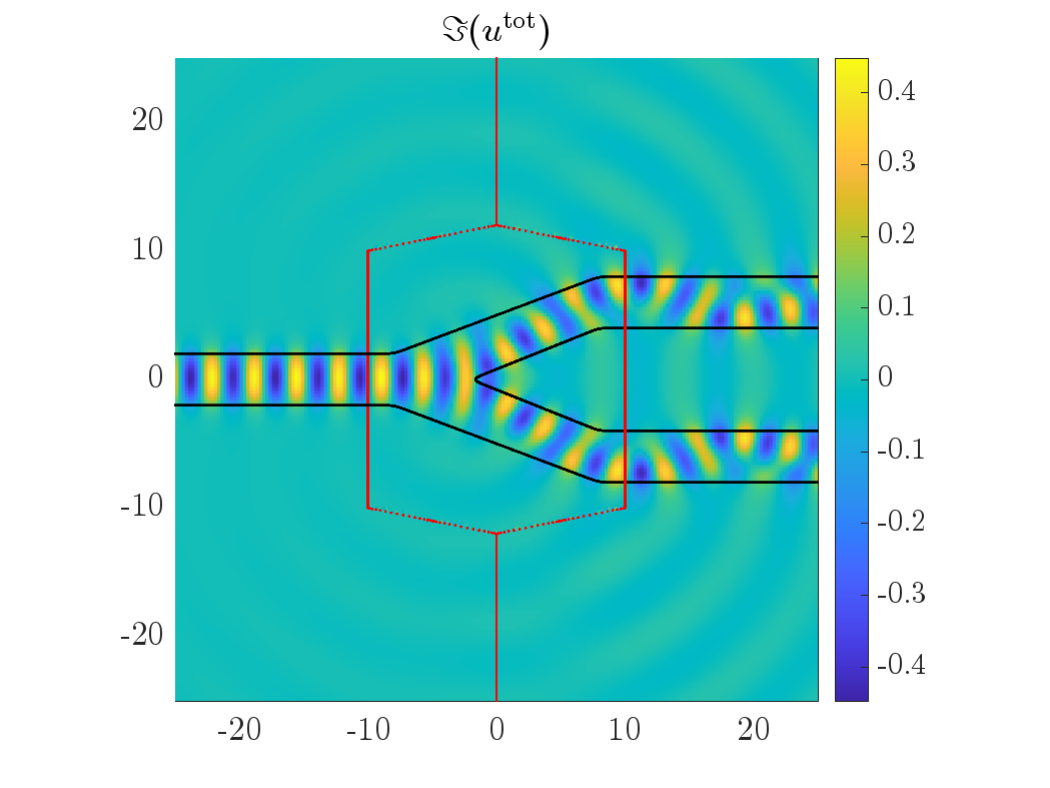}
    \caption{Simulation of a wave-guide junction.}
    \label{fig:egg_soln}
\end{figure}

In this section, we describe how to simulate a junction between more than 2 wave-guides. We model the case where a single wave-guide is split into two parallel wave-guides (see~\figref{fig:egg_soln}). In this case, we split the domain into three regions: the left and right wave-guides, and a compact transition region. 

For illustrative purposes, we focus on the case where all wave-guides have the same wave-number. This method could be extended to the case where wavenumber changes smoothly using standard integral equation techniques. To solve this problem, we compute the domains Green's functions for the left and center domains using the techniques described in Sections~\ref{sec:single_guide} and~\ref{sec:leaky_end}. The domain with two parallel wave-guides can be handled using the same technique, except that~$\tilde\gamma$ has four disjoint pieces. We use these Green's function to represent the solution~$u$ as
\begin{equation}
    u_{l,m,r}= \cS_{q_{l,m,r}}[\sigma] + \cD_{q_{l,m,r}}[\tau], \label{eq:egg_rep}
\end{equation}
where the integral operators are now supported on the union of the all the domain boundaries (the red lines in \figref{fig:egg_soln}). If we plug this representation into the continuity conditions along those boundaries, then we  arrive at an integral equation that is equivalent to~\eqref{eqn2.13.001}. In this case, however, the operators~$\cD_{q_{l,m,r}}$ and~$\cS_{q_{l,m,r}}'$ are no longer zero on the boundary of the compact region.

We discretize this integral equation using similar discretization techniques to those described, and  solve the system to compute the field generated by an incoming wave-guide mode. The resulting field is shown in~\figref{fig:egg_soln}. An analytic solution tests indicates that the solution is accurate to~$6$ digits everywhere away from the interfaces and the method described in Section~\ref{sec:NRG_calc} shows that the solution satisfies~\eqref{eq:power_bal} to 5 digits. 

\begin{figure}
    \centering
    \includegraphics[width=12cm]{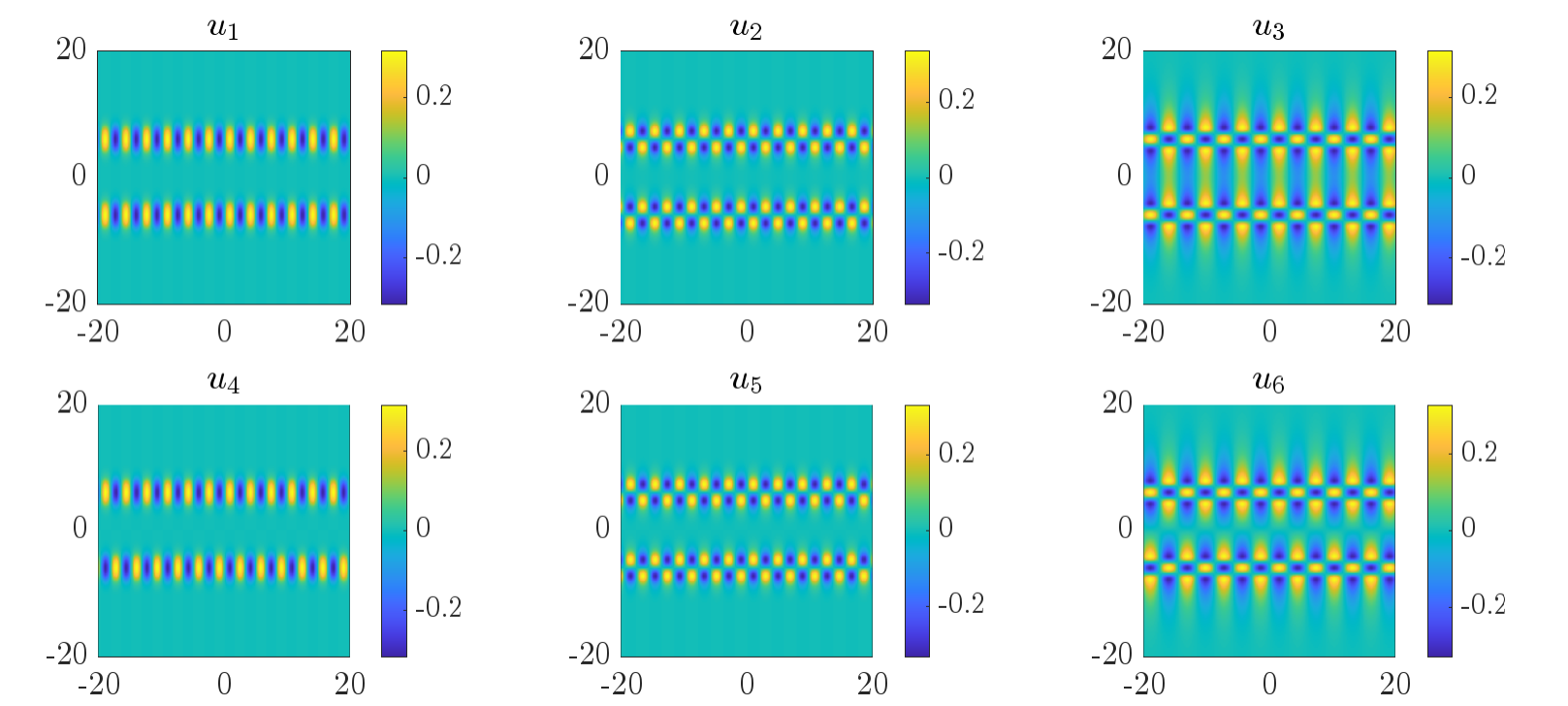}
    \caption{The wave-guide modes for the two parallel wave-guides used in~\figref{fig:egg_soln}.}
    \label{fig:two_modes}
\end{figure}

\begin{figure}
    \centering
    \includegraphics[width= 8cm]{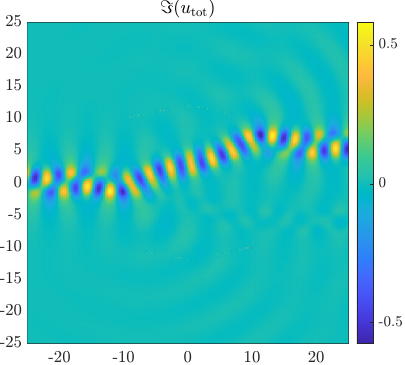}
    \caption{The combination of incoming wave-guide modes that sends the greatest power to the upper wave-guide.}
    \label{fig:power_up}
\end{figure}

With this solver, we can tune the input signal to maximize the power sent down either channel. To do this, we must find the wave-guide modes for a pair of parallel channels. It is not hard to show that these wave-guide modes are  perturbations of the modes for a single wave-guide. Indeed, if the separation of the wave-guides is~$2h$, then the wave-guide frequencies are the roots of
\begin{multline}
    f_{e}(\xi) = \left[ k_{1;l,r}^2-\xi^2+\lp k^2-\xi^2\rp\tanh\lp h\sqrt{k^2-\xi^2}\rp\right]  \sin\lp 2d\sqrt{k_{1;l,r}^2 - \xi^2} \rp \\
    - \lp1+\tanh\lp h\sqrt{k^2-\xi^2}\rp\rp\sqrt{\xi^2-k^2}\sqrt{k_{1;l,r}^2 - \xi^2}\cos\lp2d\sqrt{k_{1;l,r}^2 - \xi^2} \rp
\end{multline}
if~$v_{r,n}$ is even or
\begin{multline}
    f_{o}(\xi) = \left[ k_{1;l,r}^2-\xi^2+\lp k^2-\xi^2\rp\coth\lp h\sqrt{ k^2-\xi^2}\rp\right]  \sin\lp 2d\sqrt{k_{1;l,r}^2 - \xi^2} \rp \\
    - \lp1+\coth\lp h\sqrt{ k^2-\xi^2}\rp\rp\sqrt{\xi^2-k^2}\sqrt{k_{1;l,r}^2 - \xi^2}\cos\lp2d\sqrt{k_{1;l,r}^2 - \xi^2} \rp,
\end{multline}
if it is odd. The functions~$f_e$ and~$f_o$ are both  perturbations of~\eqref{eq:smth_mode}. The resulting modes are shown in \figref{fig:two_modes}. With these modes in hand, we can build the scattering matrix using the method from Section~\ref{sec:scatter_mat}. Using a singular value decomposition of the scattering matrix, we pick~$u^{\In}$ to be the combination of modes that maximizes the power in the upper channel. The resulting solution is shown in \figref{fig:power_up}. Using this combination, 98.32\% of the incoming power is sent into the upper channel, 0.27\% is sent to the lower channel, 0.015\% is reflected, and 1.40\% is radiated, which shows that energy is conserved to 4 significant digits.


\subsection{Dirichlet Wave-guides}
As a final example, we consider the case where Dirichlet boundary conditions are imposed on the edges of the wave-guides. To compute the Green's function for a Dirichlet wave-guide, we use different  methods, depending on the location of~$\bx$ and~$\by$. First, if~$\bx$ and~$\by$ are not in the same section of the wave-guide, the Green's function is zero. Secondly, if~$\bx$ and~$\by$ are both above or both below the wave-guide, then Green's function can be computed by placing an image charge at the reflection of~$\by$ in the closer edge of the wave-guide. For the final case, where~$\bx$ and~$\by$ are inside the wave-guide, we compute the Green's function by solving an integral equation on the boundary of the wave-guide. Specifically, we suppose that
\begin{equation}
    G_q^{\out} = \cD_k[\rho],
\end{equation}
where the integral operator is supported on the boundary of the wave-guide. Plugging this into the boundary condition~$G_q^{\In}+G_q^{\out}=0$ on that boundary gives a second kind integral equation for~$\rho$. We solve this integral equation using the same coordinate complexification and discretization techniques described above.

We can use this Green's function to simulate the case where the Dirichlet wave-guide ends without a cap. To make this simulation, we choose the left subdomain to be a Dirichlet wave-guide and the right subdomain to be the Helmholtz equation with Dirichlet boundary conditions imposed on two line segments that form the end of the pipe. By introducing these line segments, we move the corner singularity at the ends of the pipe away from the~$x_2$-axis, which ensures that the densities in the equivalent of~\eqref{eqn2.13.001} are piecewise smooth. 

To find the Green's function for the right system, we must solve the so-called open arc problem. These problems are interesting in their own right and a lot work has been done towards their accurate solution, see the references in~\cite{Jiang2004}. In this paper, we build a simple solver for illustrative purposes. We describe our algorithm in the case the case that Dirichlet boundary condition are applied on the line segment~$[-1,1]\times\{0\}$. More line segments can be added easily.

To find the Green's function for this scattering problem, we use the representation
\begin{equation}
    G_{q_r}^{\out} = \cS_k[\rho],
\end{equation}
where~$\cS_k$ is an integral operator supported on the line segment. Plugging this representation into the Dirichlet boundary conditions leads to the first kind integral equation
\begin{equation}
    \cS_k[\rho] = -G_{q_r}^{\In}.\label{eq:open_IE}
\end{equation}
In~\cite{Jiang2004} it was observed that the solution to~\eqref{eq:open_IE} has an inverse square root singularity at each end of the line segment. We therefore introduce the scaled density
\begin{equation}
    \tilde\rho(t) =\rho(t) \sqrt{1-t^2},
\end{equation}
which can be expected to be smooth up to $t=\pm 1$. If we plug this definition into~\eqref{eq:open_IE} we see that~$\tilde\rho$ solves 
\begin{equation}
    \int_{-1}^1 G_k((t,0);(t',0)) \frac{\tilde\rho(t')}{\sqrt{1-(t')^2}} dt' = -G_{q_r}^{\In}((t,0)). \label{eq:open_IE_change}
\end{equation}
The singularity in~\eqref{eq:open_IE_change} can be removed with the change of variables~$t=\cos u$:
\begin{equation}
    \int_{0}^\pi G_k((\cos u,0);(\cos u',0)) \tilde\rho(\cos u') du' = -G_{q_r}^{\In}((\cos u',0)), \label{eq:open_IE_nice}
\end{equation}
Since this change of variables has removed the singularity, this equation can be discretized with a few equally sized Gauss-Legendre panels.

As a first kind equation, \eqref{eq:open_IE_nice} is poorly conditioned when it is finely discretized. For illustrative purposes, however, we need only a fairly coarse discretization, and so the condition number is bounded by~$4\times10^{3}$. 

\begin{figure}
    \centering
    \includegraphics[width= 15cm]{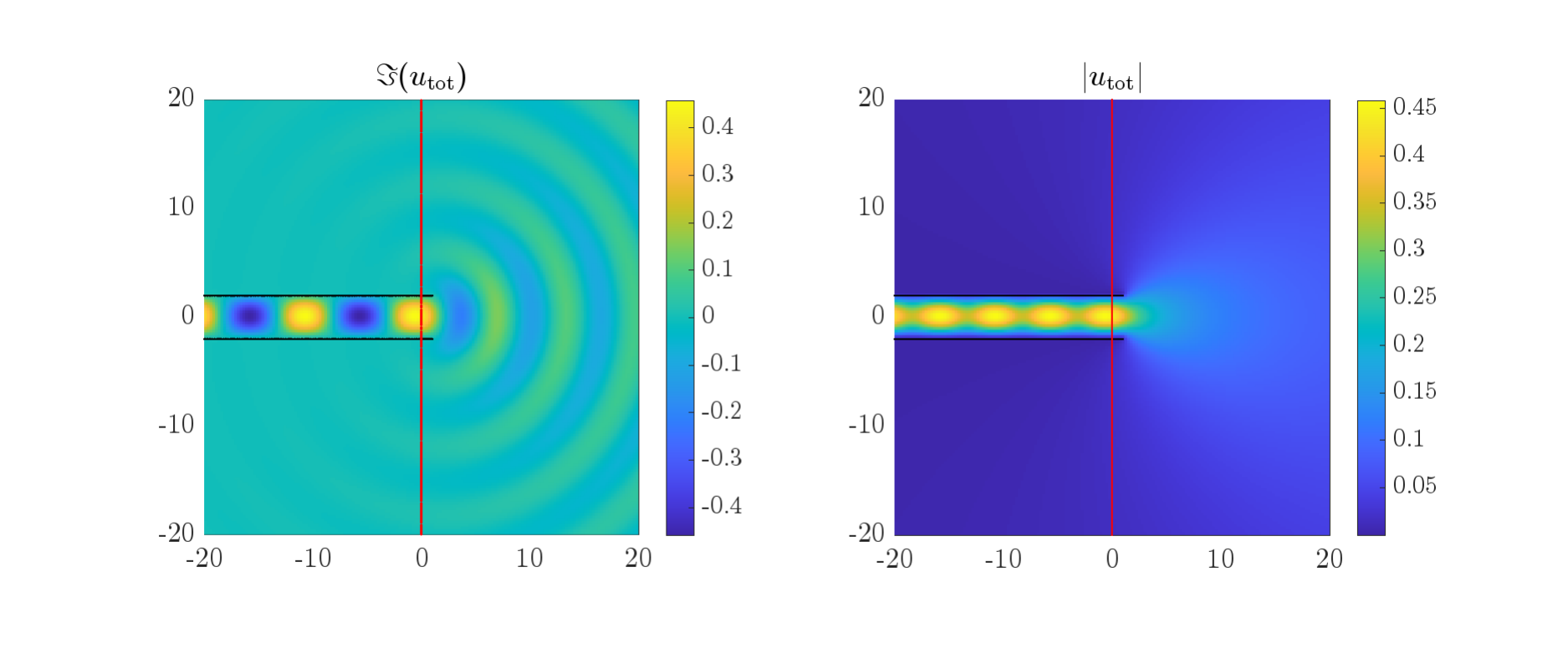}
    \caption{This figures show the imaginary part and the modulus of the solution of a terminate Dirichlet wave-guide.}
    \label{fig:Open_Dirichlet_wave-guide}
\end{figure}

By using the  methods described above to compute~$G_{q_l}$ and~$G_{q_r}$ we can reduce to solving an integral equation on the line dividing the two subdomains. \figref{fig:Open_Dirichlet_wave-guide} shows the resulting solution generated by a point source at~$(-21,0)$. An analytic solution test indicates that the solution is accurate to 6 digits everywhere away from~$\tilde\gamma$.

\begin{figure}
    \centering
    \includegraphics[width= 15cm]{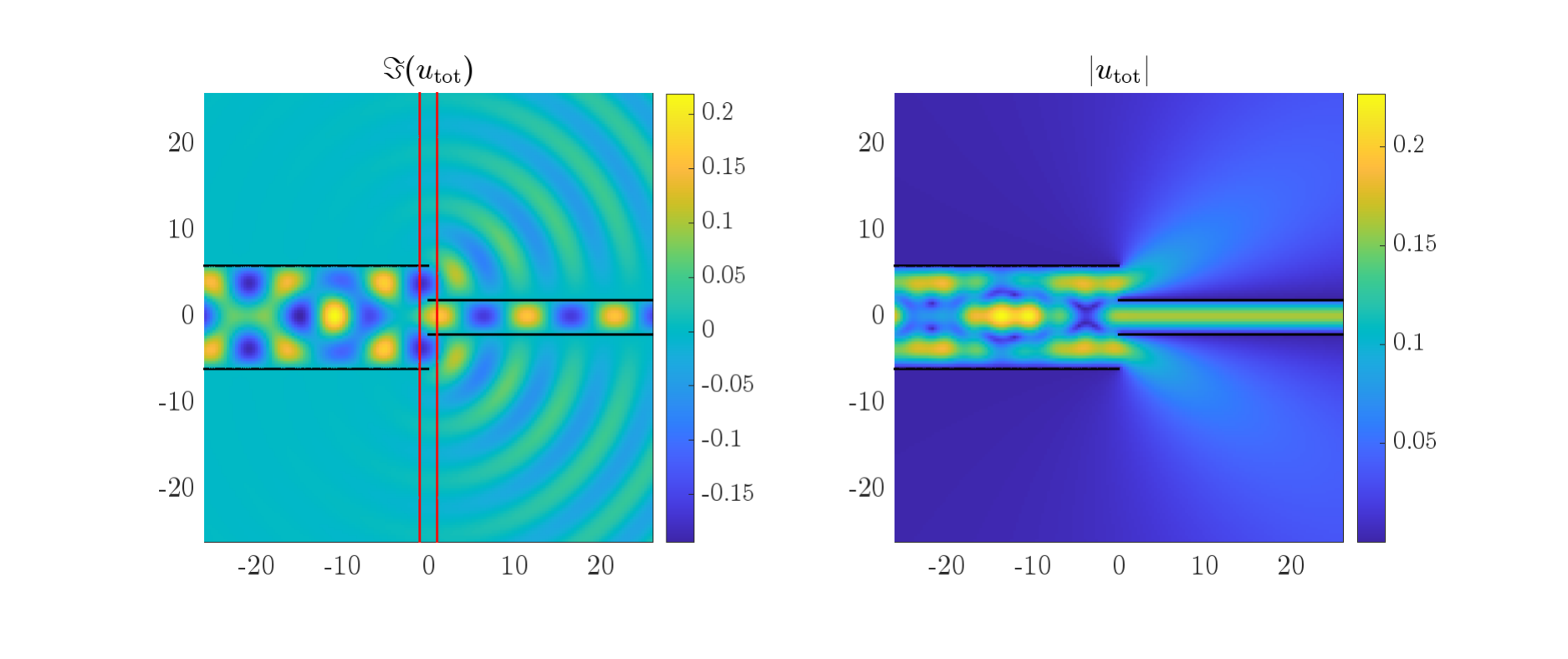}
    \caption{This figures show the imaginary part and the modulus of the solution of two Dirichlet wave-guides placed end to end.}
    \label{fig:two_Dirichlet_wave-guide}
\end{figure}

We can also use this technique to simulate the transmission of a wave between two nearby Dirichlet wave-guides. In this case, we split the problem into three regions. On either side, we have a Dirichlet wave-guide. In the middle region, we apply Dirichlet boundary conditions on a few line intervals. To solve the problem, we impose continuity conditions on two lines parallel to the~$x_2$-axis and solve the resulting system of integral equations. The resulting solution generated by a point source at~$(-31,0)$ is shown in~\figref{fig:two_Dirichlet_wave-guide}. An analytic solution test indicates that the solution is accurate to 5 digits everywhere away from~$\tilde\gamma$.

\section{Conclusion}
In this paper we have introduced the class of decomposable scattering problems on unbounded domains. We begin with an operator $L.$ A problem is decomposable if the computational domain can be split up into pieces, such that, in each unbounded component,  the operator $L$ has an outgoing fundamental solution that can be computed efficiently. We show how to use the method of fundamental solutions to reduce the problem to systems of integral equations on the union of the domain boundaries. The method becomes practicable if the domain fundamental solutions have analytic continuations that allow the resulting integral equation to be solved using coordinate complexification and standard integral equations techniques.

When implemented naively, these integral equation methods can be slow, as the evaluations of the kernel functions of the domain fundamental solutions require solving a PDE. We therefore introduce a method for building compressed representations of these domain fundamental solutions, which allow us to solve these integral equations efficiently.

This method is applied to the case of 2 semi-infinite wave-guides meeting along a common perpendicular line. We conclude with a number of other examples of decomposable problems, including a wave-guide with a  Y-junction, a wave-guide with a bend, and Dirichlet wave-guides. In future works, we plan to provide a rigorous justification for our use of coordinate complexification to solve~\eqref{eqn2.13.001} and to extend this method to other decomposable problems. 

\appendix
\section{Energy Conservation for Open Wave-guides}\label{App1}

In this appendix we give the proof of Theorem~\ref{scat_unit_thm}, which shows that the scattering matrix, at a fixed real frequency,  for an open wave-guide is a unitary map. The argument applies to any 2-dimensional wave-guide as modeled by a potential of the form given in~\eqref{eqn1}. This analysis uses concepts introduced in~\cite{EpMaSR2023}.

To begin we let $\overline{\bbR^2}$ denote the radial compactification of
$\bbR^2,$ with $\pa\overline{\bbR^2}\simeq S^1,$ the unit circle, see~\cite{EpMaSR2023}.  Suppose that
$q({\bx})$ defines an open wave-guide in the plane. This means that there is a
collection of distinct directions $\cC=\{\omega_1,\dots,\omega_N\}\subset
\pa\overline{\bbR^2}\simeq S^1$ that define the channel ends, and real valued,
compactly supported functions of 1-variable $\{q_1,\dots,q_N\},$ so that
\begin{equation}
  q({\bx})=q_0({\bx})+\sum_{j=1}^Nq_j({\bx}-\langle {\bx},\omega_j\rangle \omega_j)\varphi_+(r_j\langle {\bx},\omega_j\rangle)
\end{equation}
where $q_0\in L^{\infty}(\bbR^2)$ is compactly supported and
$\varphi_+\in\cC^{\infty}(\bbR)$ is a monotone function with
\begin{equation}
  \varphi_+(t)=
  \begin{cases}
    0\text{ for }t<1,\\
    1\text{ for }t>2.
  \end{cases}
\end{equation}
The operator $\Delta+q$ is self adjoint with domain $H^2(\bbR^2).$

We are then interested
in solutions, $u$ to
\begin{equation}
  (\Delta+q({\bx})+k^2)u({\bx})=f({\bx})\in\cS(\bbR^2),
\end{equation}
with asymptotic expansions at infinity.  By this we mean that
\begin{equation}  u(r\omega)=\frac{e^{ik r}}{r^{\frac 12}}[ a_+(\omega)+O(r^{-1})]+\frac{e^{-ik r}}{r^{\frac 12}}[ a_-(\omega)+O(r^{-1})]
\end{equation}
with
$$a_{\pm}\in\cC^{\infty}(S^1\setminus \{\omega_1,\dots,\omega_N\} )\cap
L^{\infty}(S^1).$$ In addition we assume that in a conic neighborhood, $V_j$ of
each channel end, $\omega_j\in\cC,$ in which we can write
\begin{equation}\label{eqn5}
  u({\bx})=u_c^j({\bx})+u_g^j({\bx})\text{ for }{\bx}\in V_j,
\end{equation}
where $u_g^j$ is a sum of the wave-guide modes, $\{v_{l}^j:\: l=1,\dots, M_j\},$ for the bi-infinite open wave guide
defined by $\Delta+q_j({\bx}-\langle {\bx},\omega_j\rangle \omega_j)+k^2,$ and $u_c^j({\bx})$ is the
contribution from the continuous spectrum. It satisfies
\begin{equation}\label{eqnA.6}
  u_c(r\omega)=O(r^{-\frac 12}),
\end{equation}
uniformly for $r\omega\in V_j.$ The estimate in~\eqref{eqnA.6} follows by combining the estimates in (95), (209) and (280) from~\cite{EpWG2023_2}, using the uniformity of the estimate in (280) as $\theta_{\pm}\to 0^{\pm}$. For $l=1,\dots M_j,$
\begin{equation}
  \begin{split}
  &\lim_{x_1\to
    \infty}\int_{-\infty}^{\infty}u_c(x_1,x_2)\overline{v_l^j(x_2)}dx_2=0\text{
    and }\\
   &\lim_{x_1\to
      \infty}\int_{-\infty}^{\infty}\pa_{x_1}u_c(x_1,x_2)\overline{v_l^j(x_2)}dx_2=0,
    \end{split}
\end{equation}
where we have chosen coordinates so that $\omega_j=(1,0).$

Using results from~\cite{EpWG2023_1,EpWG2023_2} and the formula in
Appendix~\ref{App2}, these asymptotics have been established in the $2d$-case. More generally
this follows from the fact that an incoming or outgoing solution to
$(\Delta+q+k^2)u=0$ can be written in a conic neighborhood, $V_j,$ of $\omega_j$
as an integral over $\pa V_j$ with respect to the outgoing resolvent
kernels, $(\Delta+q_j+k^2+ i0)^{-1}.$ This is proved in Appendix~\ref{App2} for certain
kinds of incoming data.

\begin{figure}
  \centering \includegraphics[width= 10cm]{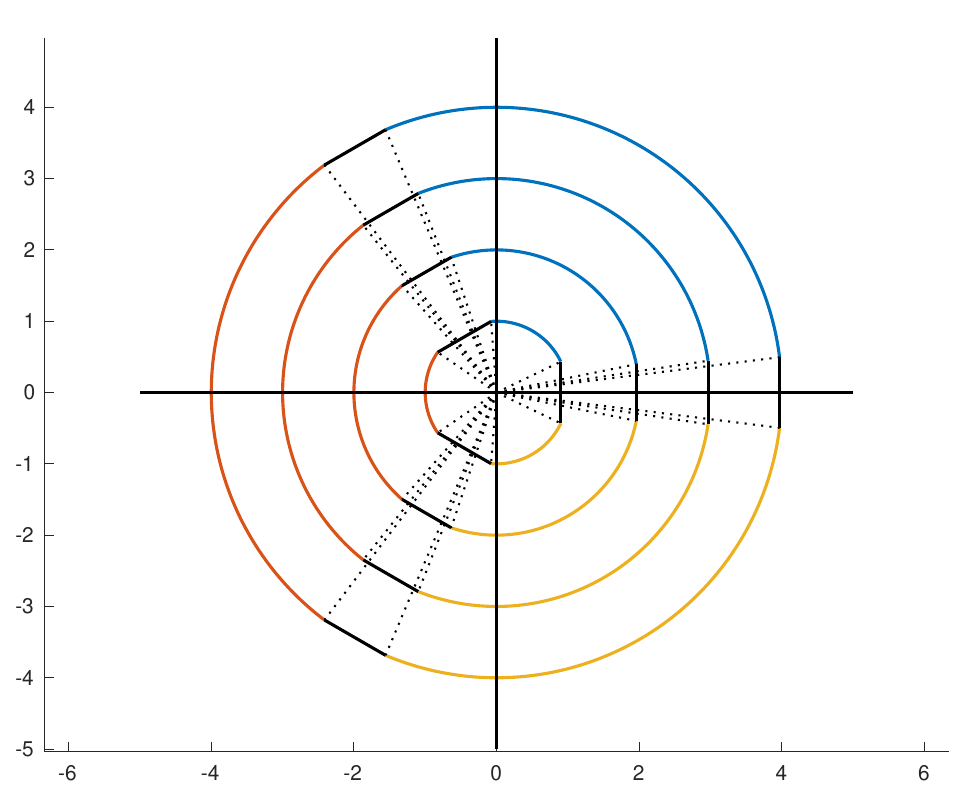}
    \caption{A collection of domains $\Omega_{R_j}$ for a wave-guide with 3 channels.}  
   \label{fig11}
\end{figure}

We now define $\{\Omega_R,:\: R>0\},$ a family of domains that exhaust $\bbR^2$
as $R\to\infty,$ adapted to the wave-guide geometry (Figure~\ref{fig11}). For an $0<\epsilon<\frac
12$ we let
\begin{equation}
  \theta_R=R^{\epsilon-1}.
\end{equation}
For $\omega_j=(\cos\theta_j,\sin\theta_j),$  let $C_{j,R}$ be the truncated cone:
\begin{multline}
  C_{j,R}=\{\lambda r(\cos(\theta_j+\theta_R),\sin(\theta_j+\theta_R))+\\
  (1-\lambda) r(\cos(\theta_j-\theta_R),\sin(\theta_j-\theta_R)):r\in
  [0,R],\,\lambda\in [0,1]\}.
\end{multline}
A point ${\bx}=r(\cos\theta,\sin\theta)\in \Omega_R$ if
\begin{equation}
  r<R\text{ and }\min\{|\theta-\theta_j|:j=1,\dots,N\}>\theta_R\text{ or }{\bx}\in
  C_{j,R}\text{ for some }j;
\end{equation}
set $\Gamma_R=\pa\Omega_R\setminus \cup_{j=1}^N\pa C_{j,R},$ and let $\pa_o
C_{j,R}=\pa C_{j,R}\cap \pa\Omega_R,$ denote the \emph{outer} boundary of
$C_{j,R}.$ Examples of $\Omega_R,$ for a range of $R,$ are shown in Figure~\ref{fig11}.

We now show that the incoming energy equals the outgoing energy using a fairly standard argument, similar to the proof of this fact given in~\cite{Melrose94}.
Suppose that we have two solutions of
\begin{equation}
  (\Delta+q+k^2)u_m=f_m\in\cS(\bbR^2),\quad m=1,2,
\end{equation}
with asymptotic expansions
\begin{equation}\label{eqn10}
      u_m(r\omega)=\frac{e^{irk}}{\sqrt{r}}a_{m+}(\omega)+\frac{e^{-irk}}{\sqrt{r}}a_{m-}(\omega)+o(r^{-\frac
        12}),\text{ for }\omega\notin\{\omega_j\}.
\end{equation}
The coefficients $\{a_{j\pm}(\omega)\}$ remain bounded as $\omega$
approaches $\omega_j$ from above and below, and that is what is understood
here. If we choose coordinates so the $\omega_j=(1,0),$ then in a conic
neighborhood, $V_j$ of $\omega_j\in\pa\overline{\bbR^2}$ we also have the expansion 
\begin{equation}\label{eqn13}
    u_m(x_1,x_2)=\sum_{l=1}^{M_j}\left[a_{ml}^{j}e^{i\xi^{j}_lx_1}v^{j}_l(x_2)+b_{ml}^{j}e^{-i\xi^{j}_lx_1}v^{j}_l(x_2)\right]
    +O(x_1^{-\frac 12}),\text{ as }\pm x_1\to \infty.
\end{equation}
The finite sum is just $\left.u_{mg}(x_1,x_2)\right|_{V_j}.$ We assume that the exponentially decaying
transverse components of the wave-guide modes,
$\{v^j_l(x_2)\},$ are normalized so that
\begin{equation}\label{eqn14}
  \int_{-\infty}^{\infty}v_l^j(x_2)\overline{v_k^j(x_2)}dx_2=\delta_{lk}.
\end{equation}

Following Melrose in~\cite{Melrose94}, we begin with the bi-linear relation
\begin{equation}
  \begin{split}
    &\lim_{R\to\infty}\int_{\Omega_R}\left[u_1\cdot
      \overline{f_2}-f_1\cdot \overline{u_2}\right]d{\bx}\\
    = &\lim_{R\to\infty}\int_{\Omega_R}\left[u_1\cdot
    \overline{(\Delta+q+k^2)u_2}-(\Delta+q+k^2)u_1\cdot \overline{u_2}\right]d{\bx}\\
 =& \lim_{R\to\infty}\int_{\pa\Omega_R}\left[u_1\cdot \overline{\pa_{\nu} u_2}-\pa_{\nu} u_1\cdot \overline{u_2}\right]ds.
  \end{split}
\end{equation}
Stokes' Theorem is used to go from the second to the third line; here $\nu$ is the outer unit normal to $\pa\Omega_R$.

First consider the  contribution of $\Gamma_R.$ For this part of
$\pa\Omega_R$ we use the expansions in~\eqref{eqn10}, which gives
\begin{equation}
     \int_{\Gamma_R}\left[u_1\cdot \overline{\pa_r u_2}-\pa_r u_1\cdot
      \overline{u_2}\right]ds
    =2ik\int_{\Gamma_R}[a_{1-}(\omega)\overline{a_{2-}(\omega)}-
      a_{1+}(\omega)\overline{a_{2+}(\omega)}]d\omega+o(1).
\end{equation}
Given our assumptions on the coefficients $a_{\pm},$ this term converges, as $R\to\infty,$ to
\begin{equation}
  2ik\int_{S^1}[a_{1-}(\omega)\overline{a_{2-}(\omega)}-
      a_{1+}(\omega)\overline{a_{2+}(\omega)}]d\omega.
\end{equation}

We now analyze the contributions of the channel ends. Fix a $j\in\{1,\dots,N\}$ and  choose coordinates
so that $\omega_j=(1,0).$ The contribution of outer portion of $\pa C_{j,R}$ takes the form
\begin{multline}
  \int_{\pa_o C_{j,R}}\left[u_1\cdot \overline{\pa_{\nu} u_2}-\pa_{\nu} u_1\cdot
    \overline{u_2}\right]ds\\
  =\int_{-R\sin\theta_R}^{R\sin\theta_R}\left[u_1\cdot \overline{\pa_{x_1} u_2}-\pa_{x_1} u_1\cdot
    \overline{u_2}\right](R\cos\theta_R,x_2)dx_2\\
  =\int_{-R\sin\theta_R}^{R\sin\theta_R}\left[u_{1g}\cdot \overline{\pa_{x_1} u_{2g}}-\pa_{x_1} u_{1g}\cdot
    \overline{u_{2g}}\right](R\cos\theta_R,x_2)dx_2+O(R^{-\frac 12}+R^{\epsilon-1}).
\end{multline}
To go from the second to the third line we use the assumptions about the
asymptotic expansion to see that, e.g.
$$u_{1c}\overline{\pa_{x_1}
  u_{2c}}(x_1,x_2)=O(x_1^{-1}),\quad u_{1g}\overline{\pa_{x_1}
  u_{2c}}(x_1,x_2)=O(x_1^{-\frac 12}e^{-\alpha |x_2|}),\text{ etc.,}$$ 
for an $\alpha>0.$ As $R\sin\theta_R\sim R^{\epsilon}$  for some $\epsilon<\frac 12,$ the estimate follows.

To compute the contribution of the remaining term we use the expansions
in~\eqref{eqn13}. The fact that $R\sin\theta_R\to\infty$ and the orthogonality
relations in~\eqref{eqn14} imply that
\begin{multline}
  \lim_{R\to\infty}\int_{-R\sin\theta_R}^{R\sin\theta_R}\left[u_{1g}\cdot \overline{\pa_{x_1} u_{2g}}-\pa_{x_1} u_{1g}\cdot
    \overline{u_{2g}}\right](R\cos\theta_R,x_2)dx_2=\\
  2i\sum_{l=1}^{M_j}
  \xi^j_l[b_{1l}^j\overline{b_{2l}^j}-
    a_{1l}^j\overline{a_{2l}^j}].
\end{multline}
Assembling the complete formula we see that
\begin{multline}
  \int_{\bbR^2}\left[u_1\cdot
      \overline{f_2}-f_1\cdot \overline{u_2}\right]d{\bx}=2ik\int_{S^1}[a_{1-}(\omega)\overline{a_{2-}(\omega)}-
    a_{1+}(\omega)\overline{a_{2+}(\omega)}]d\omega+\\
  2i\sum_{j=1}^N\sum_{l=1}^{M_j}
  \xi^j_l[b_{1l}^j\overline{b_{2l}^j}-
    a_{1l}^j\overline{a_{2l}^j}].
\end{multline}
Of particular interest is the case $u_1=u_2$ and $f_1=f_2=0;$ in this case we
see that
\begin{equation}\label{eqn21}
  k\int_{S^1}|a_{-}(\omega)|^2d\omega+
  \sum_{j=1}^N\sum_{l=1}^{M_j}
  \xi^j_l|b_{l}^j|^2=
    k\int_{S^1}|a_{+}(\omega)|^2d\omega+
  \sum_{j=1}^N\sum_{l=1}^{M_j}
  \xi^j_l|a_{l}^j|^2.
\end{equation}
This completes the proof of Theorem~\ref{scat_unit_thm}.

The scattering map in this case carries the coefficients of the incoming
asymptotics to those of the outgoing asymptotics; it is given by
\begin{multline}
  S(k):\left(\sqrt{k}a_-(\omega),\left[\sqrt{\xi_l^j}b_l^j:\: l=1,\dots,M_j,\, j=1,\dots,
    N\right]\right)\mapsto\\
  \left(\sqrt{k}a_+(\omega),\left[\sqrt{\xi_l^j}a_l^j:\: l=1,\dots,M_j,\, j=1,\dots, N\right]\right).
\end{multline}
It follows from~\eqref{eqn21} that $S(k)$ is a unitary map.

In the simplest case of two semi-infinite wave-guides meeting along common
perpendicular line we have left and
right wave-guide modes
$$\{v^l_j(x_2)e^{\pm i\xi^l_j x_1}:\: j=1,\dots, N_l\},\quad
\{v^r_j(x_2)e^{\pm i\xi^r_j x_1}:\: j=1,\dots, N_r\}.$$
For definiteness we take $\xi_j^{l,r}>0;$   normalized as before.
In this case the  unitarity statement is:
\begin{multline}
k\int_{0}^{2\pi}|a_-(\omega)|^2d\omega
+\sum_{j=1}^{N_r}\xi_j^r|b_j^r|^2+\sum_{j=1}^{N_l}\xi_j^l|b_j^l|^2=\\
k\int_{0}^{2\pi}|a_+(\omega)|^2d\omega +\sum_{j=1}^{N_r}\xi_j^r|a_j^r|^2+\sum_{j=1}^{N_l}\xi_j^l|a_j^l|^2
\end{multline}

\section{Representation of Solutions Near to Channel Ends}\label{App2}
We show that the outgoing part of a solution to $(\Delta+q+k^2)u=0$ has a
representation in a conic neighborhood of each channel end in terms of the
bi-infinite wave-guide L.A.P. resolvent $(\Delta+q_j+k^2+i0)^{-1},$ which leads
to the representations in~\eqref{eqn5}.  In this appendix we assume that the
incoming part of $u$ comes from a wave-guide mode. Similar arguments will apply
with other kinds of incoming data. We carry out this argument in $d$-dimensions.

Let $v_l^j({\bx}')e^{-i\xi^j_lx_1}$ be an incoming wave guide mode
associated to the $j$th channel, where we have chosen  coordinates so that
$\omega_j=(1,0,\dots,0).$ We choose a $\cC^{\infty}$ conic cut-off
\begin{equation}
  \Psi_{\delta}^j(r\omega)=\varphi_+(r)\psi_{\delta}(\langle\omega,\omega_j\rangle),
\end{equation}
where $\varphi_+(r)=0,$ for $r<r_j$ and $1$ for $r>2r_j,$ and
\begin{equation}
  \psi_{\delta}(t)=
  \begin{cases}
    1\text{ for }t>1-\delta,\\
        0\text{ for }t<<1-2\delta.
  \end{cases}
\end{equation}

Set $\tu_j({\bx})=\Psi_{\delta}^j({\bx})v_l^j({\bx}')e^{-i\xi^j_lx_1}, $ and observe that
\begin{equation}
  w_j=(\Delta+q+k^2)\tu_j=2\nabla\Psi_{\delta}^j\cdot\nabla(v_l^j({\bx}')e^{-i\xi^j_lx_1})+
  \Delta(\Psi_{\delta}^j)v_l^j({\bx}')e^{-i\xi^j_lx_1}\in\cS(\bbR^d).
\end{equation}
This follows as the terms $\nabla\Psi_\delta^j$ and $\Delta \Psi_\delta^j$ are supported near a cone with axis along the positive $x_1$-axis, and $v_l^j({\bx}')$ decays exponentially in this set.
This shows that $w_j$ is in the domain of the L.A.P. resolvent
$R(k^2+i0)=(\Delta+q+k^2)^{-1}$ and therefore if we let
\begin{equation}
  \tw_j=(\Delta+q+k^2)^{-1}w_j,
\end{equation}
then $\tw_j$ is outgoing, and the function
\begin{equation}
  u_{jl}^{\tot}=\tu_j({\bx})-\tw_j({\bx})
\end{equation}
satisfies $(\Delta+q+k^2)u_{jl}^{\tot}=0.$ By results in~\cite{Vasy2000}, the fact that $\tw_j$ is outgoing implies
that it has an expansion of the form
\begin{equation}\label{eqn29.1}
  \tw_j(r\omega)=\frac{e^{ikr}}{r^{\frac{d-1}{2}}}\cdot \sum_{m=0}^{\infty}\frac{a^+_m(\omega)}{r^m},
\end{equation}
where
$a^+_m\in\cC^{\infty}(\pa\overline{\bbR^d}\setminus\{\omega_j:\:j=1,\dots,N\}).$

It is also true that there is a conic neighborhood $\tV_j$ of $\omega_j$ in which
\begin{equation}
  (\Delta+q_j+k^2)\tw_j({\bx})=0.
\end{equation}
By choosing $\delta'<\delta$ we can arrange to have
$\supp\Psi^j_{\delta'}\subset\tV_j;$ as $\Psi^j_{\delta'}$ depends only on $\omega,$ for large enough $r,$ we see that
\begin{equation}
  f_j=(\Delta+q_j+k^2)[\Psi^j_{\delta'}\tw_j]=2\nabla \Psi^j_{\delta'}\cdot\nabla\tw_j+(\Delta\Psi^j_{\delta'})\tw_j=O(|{\bx}|^{-\frac{d+3}{2}}),
\end{equation}
which shows that  $f_j$ is in the domain of the $j$th L.A.P. resolvent
$(\Delta+q_j+i0)^{-1}.$

If we let $g_j=(\Delta+q_j+i0)^{-1}f_j,$ then $g_j$ is
the \emph{unique} outgoing solution to
\begin{equation}
  (\Delta+q_j+k^2)g_j=f_j.
\end{equation}
But $\Psi^j_{\delta'}\tw_j$ is also an outgoing solution to this equation, and
therefore $g_j=\Psi^j_{\delta'}\tw_j,$ which implies that
\begin{equation}
  \begin{split}
  \Psi^j_{\delta'}\tw_j({\bx})&=\lim_{R\to\infty}\int_{B_R}
  R_j(k^2+i0)({\bx};{\by})(\Delta+q_j+k^2)\Psi^j_{\delta'}\tw_j({\by})d\by\\
  &=\lim_{R\to\infty}\int_{\pa B_R}
  [R_j(k^2+i0)({\bx};{\by})\pa_r\Psi^j_{\delta'}\tw_j-\pa_r
    R_j(k^2+i0)({\bx};{\by})\Psi^j_{\delta'}\tw_j]dS_{\by}\\
 &\phantom{mmmmmm}+\Psi^j_{\delta'}\tw_j({\bx}).
  \end{split}
\end{equation}
Therefore
\begin{equation}
  \lim_{R\to\infty}\int_{\pa B_R}
  [R_j(k^2+i0)({\bx};{\by})\pa_r\Psi^j_{\delta'}\tw_j-\pa_r
    R_j(k^2+i0)({\bx};{\by})\Psi^j_{\delta'}\tw_j]dS_{\by}=0.
\end{equation}
Away from the channel ends this statement follows from~\eqref{eqn29.1}, but is non-trivial near
to $\omega_j.$ Note also that $\supp(\Psi^j_{\delta'}\tw_j)\subset\tV_j,$ which is
disjoint from $\cC\setminus\{\omega_j\}.$

With this understood we can now use Green's formula in the cone $\tV_j$ to
conclude that, for ${\bx}\in\tV_j,$
\begin{equation}
  \begin{split}
    0&=\lim_{R\to\infty}\int_{\tV_j\cap B_R}R_j(k^2+i0)({\bx};{\by})(\Delta+q_j+k^2)\tw_j({\by})d{\by}\\
    &=\lim_{R\to\infty}\Bigg[\int_{\pa\tV_j\cap
      B_R}[R_j(k^2+i0)({\bx};{\by})\pa_{\nu}\tw_j({\by})-\pa_{\nu_{\by}}R_j(k^2+i0)({\bx};{\by})\tw_j({\by})]dS_{\by}+
    \\&\int_{\tV_j\cap
      \pa
      B_R}[R_j(k^2+i0)({\bx};{\by})\pa_{\nu}\tw_j({\by})-\pa_{\nu_{\by}}R_j(k^2+i0)({\bx};{\by})\tw_j({\by})]dS_{\by}\Bigg]+\tw_j({\bx})\\
    &=\int_{\pa\tV_j}[R_j(k^2+i0)({\bx};{\by})\pa_{\nu}\tw_j({\by})-\pa_{\nu_{\by}}R_j(k^2+i0)({\bx};{\by})\tw_j({\by})]dS_{\by}+\tw_j({\bx}).
  \end{split}
\end{equation}
It is easy to see that the integral over $\pa\tV_j$ is absolutely convergent,
and above we showed that the integral over $\tV_j\cap   \pa B_R$ goes to zero as
$R\to\infty.$ This gives the desired representation formula:
\begin{equation}
 \tw_j({\bx}) =\int_{\pa\tV_j}[\pa_{\nu_{\by}}R_j(k^2+i0)({\bx};{\by})\tw_j({\by})-R_j(k^2+i0)({\bx};{\by})\pa_{\nu}\tw_j({\by})]dS_{\by},
\end{equation}
for ${\bx}\in\tV_j.$

All that is required to obtain the representation in~\eqref{eqn5} is the
description of the bi-infinite wave-guide, L.A.P. resolvent kernels,
$R_j(k^2+i0)({\bx};{\by}),$ as the sums of a continuous spectral part and the wave-guide
contribution. In the $2d$-case this is given in~\cite{EpWG2023_1,
  EpWG2023_2}. 

\bibliography{references}

\end{document}